\documentclass[final,leqno,onefignum,onetabnum,pdftex]{siamltexmm}

\usepackage{amsmath,amssymb}

\newcommand{\peq}{p_{{\rm eq}}}

\newcommand{\comment}[1]{}

\title{Nonparametric Uncertainty Quantification for Stochastic Gradient Flows\footnotemark[4]}
\author{Tyrus Berry\footnotemark[2]\ \and John Harlim\footnotemark[2]\ \footnotemark[3]\ \footnotemark[5]}
\date{\today}

\begin{document}

\maketitle
\newcommand{\slugmaster}{%
\slugger{juq}{xxxx}{xx}{x}{x--x}}

\footnotetext[2]{Department of Mathematics, the Pennsylvania State University, University Park, PA}
\footnotetext[3]{Department of Meteorology, the Pennsylvania State University, University Park, PA}
\footnotetext[4]{This work is partially supported by the Office of Naval Research Grants MURI N00014-12-1-0912}
\footnotetext[5]{The research of JH is also partially supported by the Office of Naval Research Grants N00014-13-1-0797 and the National Science Foundation DMS-1317919.}

\begin{abstract}
This paper presents a nonparametric statistical modeling method for quantifying uncertainty in stochastic gradient systems with isotropic diffusion. The central idea is to apply the diffusion maps algorithm to a training data set to produce a stochastic matrix whose generator is a discrete approximation to the backward Kolmogorov operator of the underlying dynamics. The eigenvectors of this stochastic matrix, which we will refer to as the diffusion coordinates, are discrete approximations to the eigenfunctions of the Kolmogorov operator and form an orthonormal basis for functions defined on the data set.  Using this basis, we consider the projection of three uncertainty quantification (UQ) problems (prediction, filtering, and response) into the diffusion coordinates. In these coordinates, the nonlinear prediction and response problems reduce to solving systems of infinite-dimensional linear ordinary differential equations.  Similarly, the continuous-time nonlinear filtering problem reduces to solving a system of infinite-dimensional linear stochastic differential equations. Solving the UQ problems then reduces to solving the corresponding truncated linear systems in finitely many diffusion coordinates.  By solving these systems we give a model-free algorithm for UQ on gradient flow systems with isotropic diffusion. We numerically verify these algorithms on a 1-dimensional linear gradient flow system where the analytic solutions of the UQ problems are known.  We also apply the algorithm to a chaotically forced nonlinear gradient flow system which is known to be well approximated as a stochastically forced gradient flow.
\end{abstract}




\section{Introduction}
An important emerging scientific discipline is to quantify the evolution of low-order statistics of dynamical systems in the presence of uncertainties due to errors in modeling, numerical approximations, parameters, and initial and boundary conditions \cite{smith:13}. While many uncertainty quantification (UQ) methods have been proposed, they rely on some knowledge about the underlying dynamics, at least at the coarse grained levels. One class of popular UQ methods is low-order truncation modeling based on projection of the dynamics onto some basis \cite{smith:13,mk:10}, such as the principal component analysis \cite{jolliffe:86,kutz:13}, the truncated polynomial chaos basis \cite{najm:09,mk:10}, the dynamically orthogonal basis \cite{sl:09}, the Nonlinear Laplacian Spectral Analysis (NLSA) basis \cite{gm:12}, and recently, the ``ROMQG" which stands for Reduced Order Modified Quasilinear Gaussian" method  \cite{sm:13} which was carefully designed for predicting statistics of turbulent systems. 

This paper considers UQ problems in which one has no explicit parametric form for the underlying dynamics. Given only a time series of the underlying dynamics, our goal is to devise UQ methods based on nonparametric modeling, applying ideas from diffusion maps \cite{BN,diffusion,VB}. While the classical Cauchy problem for solving the backward Kolmogorov equation is to find the semigroup solutions of this PDE, the diffusion maps technique can be interpreted as the inverse of this Cauchy problem. Namely, given a realization of the associated stochastic process, we construct a stochastic matrix (also known as Markov matrix), whose generator is a discrete approximation to the backward Kolmogorov operator. Numerically, the stochastic matrix is constructed by evaluating an appropriate kernel function on all pairs of data points and carefully re-normalizing the resulting matrix to account for sampling bias in the data set and to insure the Markov property. In this paper, we apply the recently developed variable bandwidth diffusion kernel \cite{VB}, which generalizes the original result of diffusion maps introduced in \cite{diffusion} to non-compact manifolds. Our choice is motivated by the fact that this variable bandwidth kernel is more robust and it significantly improves the operator estimation in areas of sparse sampling, e.g. on the tail of a distribution \cite{VB}.  The UQ methods discussed in this paper are only applicable for a class of nonlinear dynamical systems which can be described by stochastically forced gradient flows, since this is the class of problems for which one can learn the (backward) Kolmogorov operator from the diffusion maps algorithm \cite{diffusion,VB}.  In the conclusion we discuss the possibility of applying the UQ framework developed here to more general systems.

The representation of the probabilistic forecasting problem in diffusion coordinates was first noted in \cite{gradFlow1,gradFlow2}, however this forecasting approach was never demonstrated.  In fact, such forecasting would have been difficult for the non-compact examples which we will consider, because they used fixed bandwidth kernels in \cite{gradFlow1,gradFlow2} which would not recover the true Kolmogorov operator, as shown in \cite{VB}. In fact, the goal of \cite{gradFlow1,gradFlow2} was not to solve the forecasting problem, but instead to a low-dimensional representations of complex dynamical systems.  We should note that a `locally-scaled' version of the approach in \cite{gradFlow1,gradFlow2} was taken in \cite{gradFlow3}, which used an ad-hoc variable bandwidth function.  While the ad-hoc bandwidth of \cite{gradFlow3} is valid for finding low-dimensional coordinates, it would not be useful for the forecasting problem since, as shown in \cite{VB}, the bandwidth function would affect the estimated operator.  The previous work in \cite{gradFlow1,gradFlow2,gradFlow3} all focused on using diffusion coordinates as a nonlinear map to obtain a low-dimensional representation of the observed dynamics. In this paper we take a significantly different perspective, namely, we treat the diffusion coordinates as a basis for smooth functions and we represent probability densities in this basis for the purposes of UQ.

We will consider three UQ problems, namely the problems of forecasting, filtering \cite{mh:12}, and response \cite{mag:05}. While most stochastic projection methods chose the projection basis coordinate \cite{smith:13,mk:10} depending on the nature of the problems (such as the geometry of the state space or the distribution of the parameters), here we consider solving these UQ problems in the diffusion coordinates. For diffusion processes on nonlinear manifolds, this data-driven basis is the most natural choice since these eigenfunctions implicitly describe both the dynamics and geometry of the underlying dynamical system. For example, in this coordinate basis, the forecasting problem reduces to solving a diagonal infinite-dimensional systems of linear ODEs, as originally noted in \cite{gradFlow1}. For continuous-time nonlinear filtering problems, the solutions are characterized by the unnormalized conditional distribution which solves the Zakai equation \cite{bc:09}. In diffusion coordinates, we shall see that this nonlinear filtering problem reduces to solving a system of infinite-dimensional multiplicatively forced linear stochastic differential equations. Finally, we consider the response problem studied in \cite{mag:05}, which applied the fluctuation-dissipation theorem to estimate the statistical linear response of perturbed systems given the data set of the unperturbed system and the functional form of the external perturbation. Finding the response in our framework requires the diffusion coordinates of the perturbed system.  By assuming the perturbation is given by a known change in the potential function, we show that this perturbed basis can be constructed using an appropriate kernel, evaluated on unperturbed data. In this basis, the corresponding nonlinear response problem reduces to solving an infinite dimensional linear system of ODEs.  In each case described above, the nonparametric UQ algorithms consist of solving truncated systems of ODEs and SDEs with finitely many diffusion coordinates.

This paper will be organized as follows: In Section~\ref{DM}, we briefly review the diffusion maps and several computational manipulations for obtaining the appropriate nonparametric probabilistic models by projecting the operators to a basis of the data set itself. In Section~\ref{UQ}, we formulate two UQ problems (prediction and filtering) in the diffusion coordinates. In Section~\ref{sec4}, we will formulate the third UQ problem (nonlinear response) in the diffusion coordinates. In Section~\ref{sec5}, we verify our numerical methods on linear and nonlinear examples. We conclude the paper with a short summary and discussion in Section~\ref{sec6}. We accompany this paper with several movies in the electronic supplementary materials, depicting the evolution of the estimated time-dependent distributions.


\section{Diffusion Kernels as Nonparametric Models}\label{DM}
Consider a state variable $x \in \mathcal{M}$ evolving on a $d$-dimensional manifold $\mathcal{M} \subset \mathbb{R}^n$ according to the stochastic gradient flow,
\begin{align}\label{SDE} dx = -\nabla U(x)dt + \sqrt{2D}dW_t, \end{align}
where $U(x)$ denotes the potential of the vector field at $x\in\cal{M}$ and $D$ is a positive scalar that characterizes the amplitude of a $d-$dimensional white noise $dW_t$ on the manifold $\cal M$. Given data $\{x_i\}_{i=1}^N$ sampled independently from the invariant measure of \eqref{SDE}, $\peq(x) \propto \exp(-U/D)$, our goal is to approximate the generator, \begin{align}\label{BKop} \mathcal{L} = D\Delta  - \nabla U \cdot \nabla , \end{align}
with a stochastic matrix constructed from the data set. In order to estimate $\peq$ from the data set, we will assume that \eqref{SDE} is ergodic, and for estimating the coefficient $D$ we will require the system to be wide-sense stationary as well.  The Markov matrix that we will construct is a nonparametric model in the sense that the structure of the manifold $\mathcal{M}$ and the form of the potential function are not assumed to be known.  In \eqref{BKop}, $\Delta$ is the Laplacian (with negative eigenvalues), and $\nabla$ is the gradient operator on $\cal M$ with respect to the Riemannian metric inherited from the ambient space $\mathbb{R}^n$.

The diffusion maps algorithm \cite{diffusion} generates a stochastic matrix by evaluating a fixed bandwidth isotropic kernel function on all pairs of data points and then renormalizing the matrix and extracting the generator.  When the data set lies on a compact manifold $\cal M$, the resulting matrix converges to $D^{-1}\mathcal{L}$ in the limit of large data for any smooth potential $U$. Recently, this result was extended by both authors to non-compact manifolds by using a variable bandwidth kernel \cite{VB}. They showed that the variable bandwidth kernels produce significantly improved estimation in areas of sparse sampling (e.g. such as rare events which occur on the tail of a distribution) and are less dependent on the choice of bandwidth. As mentioned in the introduction, we will consider the variable bandwidth kernels \cite{VB} as the key ingredient in constructing the nonparametric models in this paper.

In Section \ref{findingBK} below, we briefly review the methodology in \cite{VB} for approximating the generator $D^{-1}\mathcal{L}$. Since our goal is to approximate $\mathcal{L}$, then we need to determine $D$ from the data. In Section~\ref{findingD} will provide one method to determine $D$ using the correlation time of the data and we also mention several other methods to determine $D$.

\subsection{Approximating the Generator of Stochastic Gradient Flows}\label{findingBK}

The key to our approach is the intuition that continuous notions such as functions and operators have discrete representation in the basis of the data set itself.  Given a data set $\{x_i\}_{i=1}^N\subset \mathcal{M}\subset \mathbb{R}^n$, sampled independently from the invariant measure $\peq(x)$, a function $f$ is represented by a vector $\vec f = (f(x_1),f(x_2),...,f(x_N))^\top$.  Similarly, an integral operator $G_{\epsilon} f(x) \equiv \int_{\mathcal{M}} K(x,y)f(y)\peq(y)dV(y)$, where $dV(y)$ denotes the volume form on $\cal M$, is represented by a matrix-vector multiplication between the $N\times N$ matrix, $K_{ij} = K(x_i,x_j)$ and the $N$-dimensional vector, $\vec f$. With these definitions, the matrix product $K\vec f$ yields a vector of length $N$ with $i$-th component,
\begin{align}\frac{1}{N}\left(K \vec f \right)_i = \frac{1}{N}\sum_{j=1}^N K_{ij}\vec f_j = \frac{1}{N}\sum_{j=1}^N K(x_i,x_j)f(x_j)  \stackrel{N\rightarrow\infty}{\longrightarrow}  G_\epsilon f(x_i), \nonumber \end{align}
where the limit follows from interpreting the summation as a Monte-Carlo integral.  Thus, the matrix $K$ takes functions defined on $\{x_i\}$ to functions defined on $\{x_i\}$, so in this sense we think of $K$ as an operator written in the basis of delta functions $\{\delta_{x_i}\}$ on the data set.  Notice the crucial fact that the data is not uniformly distributed on $\cal M$, so the operator is biased by the sampling measure $\peq(y)$.  This same bias applies to inner products, if $f,g$ are functions on $\cal M$ and $\vec f, \vec g$ are their representations as vectors evaluated at $\{x_i\}$, then the dot product has the following interpretation,
\begin{align} \frac{1}{N}\vec f \cdot \vec g = \frac{1}{N} \sum_{i=1}^N f(x_i)g(x_i) \stackrel{N\rightarrow\infty}{\longrightarrow} \int_{\mathcal{M}}f(y)g(y)\peq(y)dV(y) \equiv \langle f,g\rangle_{L^2(\mathcal{M},\peq)}. \nonumber \end{align}
The previous formula shows that inner products weighted by the sampling measure $\peq$ will be easy to compute in this framework.  

With the above intuition in mind, we construct a matrix $L_{\epsilon}$ that will converge to the generator $D^{-1}\mathcal{L}$, where $\cal L$ is defined in \eqref{BKop}, in the limit as $N\to\infty$ in the sense that,
\[\lim_{N\to\infty} \sum_{j=1}^N (L_{\epsilon})_{ij}\vec f_j = D^{-1}\mathcal{L}f(x_i) + \mathcal{O}(\epsilon).\]
The theory developed in \cite{VB} shows that such an $L_\epsilon$ can be constructed with a variable bandwidth kernel,
\begin{align}
K_{\epsilon}(x,y) &= \exp\left\{\frac{-\|x-y\|^2}{4\epsilon \rho(x)\rho(y)}\right\},\label{VBkernel}
\end{align}
where the bandwidth function, $\rho$, is chosen to be inversely proportional to the sampling density, $\rho \approx \peq^{-1/2}$. This choice enables us to control the error bounds in the area of sparse sample (see \cite{VB} for the detailed error estimates). Using this kernel requires us to first estimate the sampling density $\peq$ in order to define the bandwidth function $\rho$.  There are many methods of estimating $\peq$, such as kernel density estimation methods \cite{RosenblattFBK,ParzenFBK,ScottVBK,ScottVBK2}. In this paper, we estimate $\peq$ using a variable bandwidth kernel method based on the distance to the nearest neighbors (see Section~4 of \cite{VB} for details). 

Given the kernel in \eqref{VBkernel}, we apply the following normalization (first introduced in \cite{diffusion}) to remove the sampling bias, 
\begin{align}
K_{\epsilon,\alpha}(x_i,x_j) &= \frac{K_{\epsilon}(x_i,x_j)}{q_{\epsilon}(x_i)^{\alpha}q_{\epsilon}(x_j)^{\alpha}} \mbox{, \quad where}\quad q_{\epsilon}(x_i) = \sum_{j=1}^N \frac{K_{\epsilon}(x_i,x_j)}{\rho(x_i)^d}. \nonumber
\end{align}
Note that $q_{\epsilon}$ is a kernel density estimate of the sampling distribution $\peq$ based on the kernel $K_{\epsilon}$ and this is the estimate which will be used in all our numerical examples below. Throughout this paper we will use the normalization with $\alpha=-d/4$ as suggested by the theory of \cite{VB}, where $d$ is the intrinsic dimension of the manifold $\mathcal{M}$.  We note that $d$ and $\epsilon$ can be estimated from data using the method described in \cite{VB,BGH14}. The next step is to construct a Markov matrix from $K_{\epsilon,\alpha}$ by defining,
\begin{align}
\hat K_{\epsilon,\alpha}(x_i,x_j) = \frac{K_{\epsilon,\alpha}(x_i,x_j)}{q_{\epsilon,\alpha}(x_i)}\mbox{, \quad where}\quad q_{\epsilon,\alpha}(x_i) = \sum_{j=1}^N K_{\epsilon,\alpha}(x_i,x_j). \nonumber
\end{align}
Finally, the discrete approximation to the continuous generator $D^{-1}\cal L$ is given by,
\begin{align}
L_{\epsilon}(x_i,x_j) &= \frac{\hat K_{\epsilon,\alpha}(x_i,x_j)-\delta_{ij}}{\epsilon\rho(x_i)^2}.\label{approxgenerator}
\end{align}  

By computing the eigenvectors $\vec \varphi_i$ and eigenvalues $\lambda_i$ of $L_{\epsilon}$, for $i=0,1,\ldots$, we are approximating the eigenfunctions and eigenvalues of $D^{-1}\mathcal{L}$, respectively. Since $D$ is a scalar constant, obviously, the eigenvectors $\vec \varphi_i$ and eigenvalues $D\lambda_i$ are discrete approximations of the eigenfunctions and eigenvalues of the continuous generator $\cal L$, respectively. We should note that this trivial argument does not hold when $D$ is a general non-diagonal matrix describing anisotropic diffusion on $\cal M$ with respect to the Riemannian metric inherited from the ambient space $\mathbb{R}^n$. In that case, we suspect that one must use a different kernel function, such as the Local Kernels defined in \cite{bs:14} which we discuss in the conclusion.

In addition to the generator $\mathcal{L}$, we can also approximate the solution semi-group $e^{\epsilon\mathcal{L}}$ by the matrix,
\[ F_{\epsilon} \equiv  I_{N\times N} + \epsilon DL_{\epsilon}, \]
where $F_{\epsilon}$ has the same eigenvectors as $L_{\epsilon}$ but with eigenvalues $\xi_i = 1+\epsilon D\lambda_i$ where $\lambda_i$ are the eigenvalues of $L_{\epsilon}$.  
For arbitrary $t$ we can approximate the eigenvalues of the solution semi-group $e^{t\mathcal{L}}$ by $\xi_i^{t/\epsilon}$ which converges to $e^{D\lambda_i t}$ as $\epsilon\rightarrow 0$.

In order to insure that the eigenvectors $\vec \varphi_i$ (which approximate the eigenfunctions $\varphi_i$) form an orthonormal basis with respect to the inner product $\langle \cdot,\cdot\rangle _{L^2(\mathcal{M},\peq)}$, we note that, 
\[ \delta_{ij} = \langle \varphi_i,\varphi_j \rangle_{L^2(\mathcal{M},\peq)} = \int_{\mathcal{M}} \varphi_i(x)\varphi_j(x)\peq(x)dV(x) = \lim_{N\to\infty} \frac{1}{N}  \sum_{l=1}^N (\vec \varphi_i)_l (\vec \varphi_j)_l = \lim_{N\to\infty} \frac{ \vec \varphi_i^\top \vec \varphi_j}{N}. \]
Since $\{\vec \varphi_i\}_{i=0}^\infty$ are eigenvectors, they are already orthogonal, thus we only need to renomalize $\vec \varphi_i$ so that $\vec \varphi_i^\top \vec\varphi_i = N$.  To do this we simply replace $\vec \varphi_i$ with $\frac{\sqrt{N}}{||\vec\varphi_i||}\vec\varphi_i$. To simplify the notation below, we define $\langle \cdot,\cdot\rangle _{\peq}\equiv\langle \cdot,\cdot\rangle _{L^2(\mathcal{M},\peq)}$.

\subsection{Determining the Diffusion Coefficient}\label{findingD}

The algorithm described in Section \ref{findingBK} approximates the generator $D^{-1}\mathcal{L}$, where $\cal L$ is defined in \eqref{BKop}. Since our aim is to obtain the generator $\cal L$, then we must approximate the scalar diffusion coefficient $D$. Several nonparametric methods for estimating $D$ have been proposed. For example, \cite{cve:06} considered an optimization problem based on a variational formulation which matches information from the discrete approximation of a conditional Markov chain to the generators of the forward and backward Kolmogorov operators.  Another method employed in \cite{mg:12} used the one-lag correlations to approximate the diffusion coefficient. The problem with this estimate is that when the true dynamics are not described by a gradient flow, we would like to use a gradient flow to approximate the dynamics on long time scales (see Section \ref{sec52} for an example).  If we estimate $D$ using the one-lag correlation as in \cite{mg:12}, we will only match the short-time correlation of the dynamics, but we are interested in the long-time correlation which is better captured by the correlation time.  

In this section, we introduce a method to determine $D$ using the correlation time of a one-dimensional observable $S(x(t))$, which results from applying the observation function $S:\mathcal{M}\to\mathbb{R}$ to the multidimensional time series, $x(t)$.  By computing the correlation time of this observable for the estimated dynamical system, $D^{-1}\mathcal{L}$, and comparing to the empirical correlation time estimated from the training data set, we will be able to extract the intrinsic parameter $D$. This approach generalizes the Mean Stochastic Model (MSM) of \cite{mgy:10,mh:12}. Intuitively, in the previous section we used the invariant measure to implicitly determine the potential function (through the kernel based generator), and in this section we use the correlation time to fit the stochastic forcing constant $D$.  The correlation time is defined by,
\begin{align}
T_c = \int_0^\infty C(\tau)C(0)^{-1}d\tau,\nonumber
\end{align}
where $C(\tau) \equiv \langle S(x(t+\tau))S(x(t)) \rangle$ is the correlation function of the one-dimensional time series, $S(x(t))$.  The correlation function can be determined from the data by averaging over $t$, or by the inverse Fourier transform of the power spectrum, which follows from the Wiener-Khinchin formula, $C(\tau) = \mathcal{F}^{-1}\left(\|\mathcal{F}(S(x(t)))\|^2\right)$. For small data sets we found the Wiener-Khinchin approach to be more robust and this is the approach used in the examples in Section \ref{sec5}.  We note that the observation function $S$ can be any centered functional on the data set as long as the observed process is stationary and is not identically zero (which guarantees that the correlation time is well defined). Note that in the case of anisotropic diffusion, one would need to compute the entire correlation matrix $\mathbb{E}[x(t-\tau)x(t)^\top]$ as a function of $\tau$. Since we assume the diffusion is isotropic, we will only need a single correlation statistic, and the correlation time of the time series $S(x(t))$ will be sufficient.  

Once $T_c$ is estimated from the data, we need to find the value of $D$ which makes \eqref{SDE} have correlation time $T_c$.  In order to do this, we will show that $T_c$ can be estimated from the eigenvalues and eigenfunctions of $D^{-1}\mathcal{L}$.  Let $e^{\tau\mathcal{L}}$ be the semigroup solution for the backward equation. Following \cite{mag:05} we have,
\[ C(\tau) = \langle S(x(t+\tau)) S(x(t)) \rangle  = \int_{\mathcal{M}} \big(e^{\tau\mathcal{L}}(S)(x)\big) S(x) \peq(x)dV(x). \]
Writing $S = \sum_i \langle S,\varphi_i\rangle_{\peq}\varphi_i$ in the eigenbasis of $\mathcal{L}$, we have,
\begin{align}\label{correlationfunction} C(\tau) = \int_{\cal M} \sum_i  e^{D\lambda_i\tau}\langle S,\varphi_i\rangle_{\peq}^{\top} \varphi_i(x) S(x)\peq(x) dV(x) = \sum_i e^{D\lambda_i\tau} \langle S,\varphi_i \rangle _{\peq}^2. \end{align}
Note that since $\lambda_0=0$ and $\varphi_0=1$, $\langle S,\varphi_0\rangle _{\peq} = \int_{\mathcal{M}} S(x)\peq(x)dV(x) = \mathbb{E}_{\peq}[S(x)]$ so we require $S(x)$ to be centered (otherwise the integral of the first term diverges).  Noting that $\lambda_i < 0$ for $i>0$, we can compute the correlation time analytically as,
\[ T_c = \int_0^{\infty} C(\tau)C(0)^{-1}d\tau = - \frac{\sum_{i\geq 1} (D\lambda_i)^{-1}\langle S,\varphi_i \rangle _{\peq}^2}{\sum_{j\geq 1}\langle S,\varphi_j \rangle _{\peq}^2}, \]
and therefore,
\begin{align}\label{Dformula} D = - \frac{1}{T_c}\frac{\sum_{i\geq 1} \lambda_i^{-1}\langle S,\varphi_i \rangle _{\peq}^2}{\sum_{j\geq 1}\langle S,\varphi_j \rangle _{\peq}^2}. \end{align}
This gives us an equation for the diffusion coefficient $D$ which matches any given correlation time $T_c$.  Using the empirical correlation time allows us to find the diffusion coefficient of the system defined by the data.  We should note that the correlation function, $C(\tau)$, and hence the correlation time $T_c$, are not intrinsic to the manifold $\mathcal{M}$ due to the dependence on the observation function, $S$.  However, the intrinsic diffusion constant $D$ can be recovered from the ratio \eqref{Dformula} since both formulas for $T_c$ are based the same function $S$.

For one-dimensional linear problems, such as the Ornstein-Uhlenbeck process, \eqref{Dformula} reduces to the mean stochastic model (MSM) introduced in \cite{mgy:10,mh:12}. To see this, let the potential function be $U(x) = -\alpha x^2/2\in\mathbb{R}$, so that the gradient flow system in \eqref{SDE} is a one-dimensional Ornstein-Uhlenbeck (OU) process and the eigenfunctions of $\mathcal{L}$ are the Hermite polynomials $\varphi_i(x) = H_i(x)$.  All the inner products $\langle S,\varphi_i\rangle _{\peq} = \langle x,H_i(x) \rangle _{\peq} = \langle H_1(x),H_i(x)\rangle _{\peq}$ are zero except for $i=1$, in which case the inner product is $1$ and the associated eigenvalue is $\lambda_1 = \alpha/D$ so that \eqref{Dformula} becomes $\alpha = -1/T_c$ which is the MSM formula.  For  nonlinear problems, the correlation function \eqref{correlationfunction} for the gradient flow system is more accurate than that of the MSM fit as we will show in Section \ref{sec5}.

To numerically estimate \eqref{Dformula}, we first approximate,
\[ \left<S,\varphi_i\right>_{\peq} \approx \frac{1}{N}\sum_{l=1}^N S(x_l)\varphi_i(x_l) = \frac{1}{N}S(x)^\top \vec{\varphi}_i, \]
where $S(x)$ is an $N\times 1$ matrix with $l-$th row given by, $S(x_l)$. We then approximate the diffusion coefficient as,
\[ D \approx -\frac{1}{T_c} \frac{\sum_{i=1}^M \lambda_i^{-1} (S(x)^\top \vec\varphi_i)^2}{\sum_{i=1}^M  (S(x)^\top \vec\varphi_i)^2}, \]
which estimates \eqref{Dformula} with summations over $M$ modes.  We obtained the best empirical results with the observation function $S(x(t)) \equiv \sum_{j=1}^n (x(t))_j - \mathbb{E}[(x(t))_j]$, which results from summing the (centered) coordinates of the data point $x(t)$. Of course, this observable is vulnerable to certain pathological examples, such as a two dimensional time series with $x_2 = -x_1$. A more robust choice for $S$ would be the norm of the multi-dimensional time series $x(t)$, however we found better results with the summation. 

To verify the formula for estimating $D$, we simulate the Ornstein-Uhlenbeck process with true $D=1$ to produce a time series $\{x_i\}_{i=1}^{T}$ with $T=1000000$ and discrete time step $\Delta t=0.01$.  We estimate the correlation function as a simple average of $\frac{1}{T-j}\sum_{i=1}^{T-j} S(x_{i+j})S(x_i)$ for lags $j$ starting at $j=1$ and increasing until the first value of $j$ for which the average was negative.  Using the trapezoid rule to estimate the integral over all the values of the shift, $j$, we estimate the correlation time to be $T_c \approx 1.0238$.  We then subsample every $50$-th data point to produce a time series of length $20000$, the samples of which are approximately independent samples of the invariant measure. We apply the diffusion maps algorithm with the variable bandwidth kernel to the subsampled data to estimate the operator $D^{-1}\mathcal{L}$ as well as the first $M=500$ eigenvalues and eigenfunctions.  Approximating the formula \eqref{Dformula} as described above, we found $D \approx 1.0073$.  In order to verify that this result for $D$ depends only on the intrinsic geometry, we then map the time series $\{x_i\}$ into $\mathbb{R}^3$ with the isometric embedding,
\[ x_i \mapsto \mathcal{F}(x_i) = (\sin(2x),\cos(2x),x)^\top /\sqrt{5}. \]
One can easily check that $\mathcal{F}$ is an isometry since $D_{x}\mathcal{F}^\top D_x\mathcal{F} \equiv 1$. We then repeat the above procedure for the time series $S \circ \mathcal{F}(x_i)$ and found $T_c \approx 0.7571$ and $D \approx 1.0186$.  Notice that the estimates for $D$ are very similar, whereas the estimates for $T_c$ are different. The isometric change in the geometry, $\mathcal{F}$, has decreased the correlation time for the observable $S \circ \mathcal{F}$, however the intrinsic variable $D$ is not effected.  This is because $D$ is approximated in \eqref{Dformula} as a ratio between two different methods of estimating the correlation time, and the isometry has the same effect on both estimates of the correlation time.


\section{Forecasting and Filtering in Diffusion Coordinates}\label{UQ}

In this section we formulate two UQ problems, namely forecasting and filtering, in diffusion coordinates.  These two problems involve pushing a probability measure forward in time, and therefore they involve the Fokker-Planck operator $\mathcal{L}^*$.  
By projecting these infinite dimensional systems onto the eigenfunctions of the Fokker-Planck operator, we find infinite dimensional linear systems which govern the evolution on the projected coordinates.  We can then approximate the solution by truncating the linear system to describe finitely many coordinates.

For gradient flows, it is easy to determine the eigenvalues and eigenfunctions of the Fokker-Planck operator, $\mathcal{L}^*$, from the eigenvalues and eigenfunctions of the generator $\cal L$. To see this, note that $\peq\propto e^{-U/D}$, and,
\begin{align}\label{conj}  \frac{1}{\peq}\mathcal{L}^*(f\peq) &= \frac{1}{p} \textup{div}\left(D \nabla(f\peq) + f\peq \nabla U \right) = \frac{1}{\peq}\textup{div}\left(D\peq\nabla f +D f\nabla \peq + f\peq\nabla U \right) \nonumber \\
&= \frac{1}{\peq}\textup{div}(D\peq\nabla f) = D\Delta f + D\nabla f \cdot \frac{\nabla \peq}{\peq} = D\Delta f - \nabla f \cdot \nabla U = \mathcal{L}f,
\end{align}
since $D\nabla \peq = - \peq\nabla U $. Also, it follows from \eqref{conj} that $\mathcal{L}^*f = \peq\mathcal{L}(f/\peq)$.  The two operators have the same eigenvalues and it is easy to show that the semi-group solutions are related by $e^{t\mathcal{L}^*} f= \peq e^{t\mathcal{L}}(f/\peq)$ and the eigenfunctions $\psi_i$ of $\mathcal{L}^*$ and $e^{t\mathcal{L}^*}$ are given by $\psi_i = \peq \varphi_i$ where $\varphi_i$ are the eigenfunctions of $\mathcal{L}$ and $e^{t\mathcal{L}}$. From the orthonormality of $\varphi_i$, it is easy to deduce that,
\begin{align}\langle \psi_i,\psi_j \rangle_{1/\peq} =\langle \psi_i,\varphi_j \rangle = \langle \varphi_i,\varphi_j \rangle_{\peq}  =\delta_{ij}.\label{innerproductrule}
\end{align} 
Of course, since these eigenfunctions are all represented by $N$-dimensional eigenvectors (which represent evaluation on the data set $\{x_j\}_{j=1}^N$) the actual computation of any inner product will always be realized by rewriting it with respect to $\peq$, following the rule in \eqref{innerproductrule}.

\subsection{Nonlinear forecasting}
The \emph{forecasting} problem is, given an arbitrary initial distribution $p_0(x)$ for a state variable $x$, find the density $p(x,t)$ at any future time $t$.  This will be the most straightforward problem to solve in diffusion coordinates since $p(x,t)$ solves the Fokker-Planck equation,
\begin{align}\label{forecasteq} \frac{\partial p}{\partial t} = \mathcal{L}^*p, \quad\quad p(x,0) =p_0(x).
\end{align}
Let $\psi_i$ be the eigenfunctions of $\cal L^*$ with eigenvalues $D\lambda_i$, and write $p$ in this basis as, $p(x,t) = \sum_i c_i(t) \psi_i(x)$ where $c_i = \langle p,\psi_i\rangle _{1/\peq}$.  We define the vector $\vec c = \vec c(t)$ of eigencoordinates and the diagonal matrix $\Lambda_{ii} = \lambda_i$, and writing \eqref{forecasteq} in these coordinates we have,
\begin{align}\label{forecastprojection} \frac{d\vec c}{dt} = D\Lambda \vec c, \quad\quad\vec c(0) = \langle p_0,\psi_i \rangle_{1/\peq}.
\end{align}
In order to solve the linear system in \eqref{forecastprojection}, we first need to find the initial conditions $c_i(0)$ by projecting the initial density $p(x,0)=p_0(x)$ onto $\psi_i$.  Since densities $p_0(x)$ and $\peq(x)$ are given on the training data set $x_i$, we define the discrete densities $(\vec p_0)_l = p_0(x_l)$ and $(\vec p_{{\rm eq}})_l = \peq(x_l)$.  
Following the inner product rule in \eqref{innerproductrule}, we write the projection of the initial condition as, 
\begin{align}\label{coeffest} c_i(0) =  \left<p_0,\psi_i\right>_{1/\peq} = \left<p_0/\peq,\varphi_i \right>_{\peq} \approx \frac{1}{N} \sum_{j=1}^N \frac{p_0(x_j)}{\peq(x_j)}\varphi_i(x_j), \end{align}
since this converges to the true value as $N\to\infty$. With this initial condition, each component in \eqref{forecastprojection} can be solved analytically as $c_i(t) = e^{tD\lambda_i}c_i(0)$. Numerically, we will approximate the solutions with finitely many modes, $i=0,1,\ldots,M$.  In Section \ref{QoI} we show how to use these solutions to reconstruct the density $p(x,t)$ or to estimate a quantity of interest at any time $t$.

\subsection{Nonlinear filtering}\label{nonlinearfilter}
The \emph{filtering} problem is, given a time series of noisy observations $z$ of the state variable $x$, find the posterior density $P(x,t | z(s), s\leq t)$ given all the observations $z(s)$ up to time $t$.  For an observation function $h:\mathcal{M}\rightarrow\mathbb{R}^m$, we consider continuous-time observations of the form,
\begin{align}\label{obseq} dz = h(x) \, dt + \sqrt{R} \, dW_t, \end{align}
where $dW_t$ are i.i.d Gaussian noise in $\mathbb{R}^m$.
With these observations, the evolution of the unnormalized posterior density, $p(x,t)$, defined as follows,
\[
P(x,t | z(s), s\leq t)= \frac{p(x,t)}{\int_{\cal M} p(x,t)dV(x)},
\]
 is given by the Zakai equation \cite{bc:09},
\begin{align}\label{filtereq} dp = \mathcal{L}^*p dt + ph^\top R^{-1} dz. \end{align}
Writing $p(x,t) = \sum_i c_i(t)\psi_i(x)$ the evolution becomes,
\[ \sum_i dc_i \psi_i = \sum_i D\lambda_i c_i \psi_i dt + \sum_i c_i \psi_i h^\top R^{-1} dz, \]
and projecting both sides of this equation onto $\psi_j$ we find,
\begin{align}\label{infinitesde} dc_j = D\lambda_j c_j dt + \sum_i c_i \langle \psi_i h^\top,\psi_j\rangle _{1/\peq} R^{-1} dz, \end{align}
which is an infinite-dimensional system of stochastic differential equations with multiplicative stochastic forcing. 

To numerically estimate the solution of \eqref{infinitesde}, we truncate this infinite-dimensional system of SDEs by solving only a system of $(M+1)-$dimensional SDEs for $\vec{c}=(c_0,c_1,\ldots,c_M)^{\top}$. This strategy is simply a Galerkin approximation for the Zakai equation for which we choose the diffusion coordinates as a basis rather than the Gaussian series as proposed in \cite{ar:97} or the Hermite polynomial basis as proposed in \cite{fsx:13}. We should note that for scalar OU processes, the diffusion coordinates are exactly the Hermite polynomials since these polynomials are the eigenfunctions of the backward Kolmogorov operator for the OU process.  

Defining $H_{ji} = \left<\psi_i h^\top,\varphi_j \right>$ as a $1\times m$ dimensional vector for each pair $(i,j)$; where the $k$-th component is given by $H_{ji}^k = \left<\psi_i h_k,\varphi_j \right>$, we can write the truncated $(M+1)-$dimensional system of SDEs in compact form as follows,
\begin{align}\label{filterprojection} d\vec c = D\Lambda \vec c dt +  (H\vec c) R^{-1} dz. \end{align}
To solve this system of SDEs, we first project the initial condition $p_0(x)$ to obtain $c_i(0)$ as in \eqref{coeffest} for $i=0,1,\ldots,M$. 
We then numerically solve the system of SDEs in \eqref{filterprojection} using the splitting-up method (see for example \cite{fsx:13}).  Explicitly, given $\vec c(t_{i-1})$ we compute $\vec c(t_i)$ by first using the solution to the deterministic part,
\[ \vec c_0(t_i) = \exp\left(D\Lambda \Delta t \right)\vec c(t_{i-1}), \]
and then the solution to the stochastic part,
\[ \vec c(t_i) = \exp\left(H R^{-1}dz(t_i) - \frac{1}{2}H^2 R^{-1}\vec 1 \Delta t \right)\vec c_0(t_i). \]
Here, the exponent term consists of an $M\times M$ matrix given by 
\[ HR^{-1}dz(t_i) - \frac{1}{2}H^2 R^{-1}\vec 1 \Delta t  = \sum_k H^k (R^{-1}dz(t_i))_k - \frac{1}{2}(H^k)^2 (R^{-1}\vec 1)_k \Delta t. \]
Notice that unlike the forecasting problem which had an analytic solution for any time $t$, this procedure must be iterated for each observation $dz(t_i)$.

\subsection{Quantity of Interest}\label{QoI}

In the uncertainty quantification problems described above, we reduce the PDE and SPDE problems which describe the evolution of density $p(x,t)$ into systems of ODEs in \eqref{forecastprojection} and SDEs in \eqref{filterprojection}, respectively, for the coefficients $c_i(t)$.  Using the coefficients $c_i(t)$ we can approximate the density function as follows,
\begin{align}\label{reconstruct} p(x,t) \approx p_M(x,t) =\sum_{i=0}^{M} c_i(t)\psi_i(x) =\sum_{i=0}^{M} c_i(t)\varphi_i(x)\peq(x) , \end{align}
where the coefficients $\{c_i(t)\}_{i=0}^M$ are from the solution of the $(M+1)-$dimensional systems of ODEs or SDEs.

One complication in some of the algorithms above is that the evolution of the coefficients may not preserve the normalization of the corresponding density. In fact, for the nonlinear filtering problems in \eqref{filtereq}, the solutions are unnormalized densities.  Fortunately we can use the eigenfunctions $\varphi_i$ to estimate the normalization factor of an unnormalized density from its diffusion coordinates.  Assume that $p_M$ in \eqref{reconstruct} is an unnormalized density with diffusion coordinates $c_i$ obtained from the solution of the $(M+1)-$dimensional systems of ODEs or SDEs.  We can then estimate the normalization constant as,
\[ Z = \int_{\mathcal{M}} p_M(x,t)dV(x) = \sum_{i=0}^M c_i \int_{\mathcal{M}} \varphi_i(x)\peq(x)dV(x) = \sum_{i=0}^M c_i \langle 1,\varphi_i \rangle_{\peq}. \]
where we estimate inner product, $\langle 1, \varphi_i \rangle_{\peq} \approx \frac{1}{N}\sum_{j=1}^N \varphi_i(x_j)$. Subsequently, we normalize our density by replacing $p_M$ with $p_M/Z$. In the remainder of this paper, we will refer to $p_M$ as the normalized density without ambiguity.  We also note that by reconstructing the density from the coefficients we are implicitly projecting the density into the subspace of Hilbert space spanned by the first $M+1$ eigenfunctions.  This has a smoothing effect which improves the clarity of figures and we always show reconstructed densities in the figures and videos.

In many cases, we are not interested in the entire density $p$; instead we are often interested in the expectation of an observable $A:\mathcal{M}\to\mathbb{R}$ with respect to $p$, namely $\mathbb{E}_{p(x,t)}[A(x)]$.  A common example of functionals of interest are the moments $A(x) = x^m$ for $m\in\mathbb{N}$.  Although $A$ is a functional on $\mathcal{M}$, $A$ can be given as a functional on the ambient space $\mathbb{R}^n$ because we only evaluate $A$ on the data set, which lies on the manifold $\cal M$.  Using the formula \eqref{reconstruct} for the density estimate at time $t$, we can approximate the expectation of the functional $A$ at time $t$ as,
\begin{align}\label{functionalexpectation} \mathbb{E}_{p(x,t)}[A(x)] \approx \langle A,p_M \rangle = \sum_{i=0}^{M} c_i(t) \langle A,\varphi_i\rangle _{\peq} =  \sum_{i=0}^M c_i(t) a_i = \vec c(t)^\top \vec a. 
\end{align}
The formula \eqref{functionalexpectation} shows that we can find the expectation of a functional $A$ by writing $A$ in diffusion coordinates as $a_i = \langle A,\varphi_i\rangle _{\peq}$.  In these coordinates, the expectation of $A$ is simply the inner product of the diffusion coefficients of $A$ and $p(x,t)$.  Finally, if the quantity of interest is given by a functional $A$ evaluated on the training data set as $(\vec A)_l = A(x_l)$, we approximate the inner product $a_i = \left<A,\varphi_i\right>_{\peq} \approx \frac{1}{N} \vec A^\top \vec \varphi_i$.

\section{Nonlinear Response in Diffusion Coordinates}\label{sec4}
The response problem considered in \cite{mag:05} is to determine the change in the expectation of a functional $A(x)$ as function of time, after perturbing the dynamics of a system at equilibrium.  Letting $p^\delta(x,t)$ be the density evolved according to the perturbed system, we define the response of a functional $A$ as the expected difference,
\begin{align}
\delta\mathbb{E}[A(x)](t) = \mathbb{E}_{p^\delta}[A(x)](t) - \mathbb{E}_{\peq}[A(x)],\quad t\geq 0, \label{response}
\end{align}
We assume the unperturbed system has the form,
\begin{align}
dx = F(x)\,dt + \sqrt{2D}dW.\label{unperturbed}
\end{align}
with equilibrium density $\peq$ and Fokker-Planck operator $\cal L^*$, so that $\mathcal{L}^*\peq = 0$.  We assume the perturbed system has the form,
\begin{align}
dx = (F(x)+\delta F(x,t))\,dt + \sqrt{2D}dW,\label{perturbed} 
\end{align}
with associated Fokker-Planck operator $\tilde{\mathcal{L}}^*$ so that $p^\delta$ solves the Fokker-Planck equation,
\begin{align}\label{responseeq} \frac{\partial p^\delta}{\partial t} = \tilde{\mathcal{L}}^*p^\delta = \mathcal{L}^*p^\delta + \delta\mathcal{L}^*p^\delta, \end{align}
where $\delta\mathcal{L}^*p=-\textup{div}(\delta Fp)$ denotes the Liouville operator corresponding to the external forcing with functional form $\delta F$.

This reveals that the response problem is closely related to the forecasting problem. Given the Fokker-Planck operator $\tilde{\mathcal{L}}^*$ for the perturbed system, one must solve \eqref{responseeq} with initial condition chosen to be the equilibrium distribution of the unperturbed system, $p^\delta(x,0)=\peq(x)$. Since computational cost becomes an issue, especially for higher dimensional problems, one typically approximates $p^\delta$ by solving ensemble solutions of \eqref{perturbed} with initial ensemble sampled from $\peq$, and in fact we will use this approach as a benchmark in Section \ref{dwresponsesection}. However, as in \cite{mag:05}, we are interested in cases where the Fokker-Planck operators $\mathcal{L}^*$ and $\tilde{\mathcal{L}}^*$ and also the functional form of $F$ are all unknown. We will assume that we only have access to data from the unperturbed system and the functional form of the external perturbation $\delta F$.  Of course, since the techniques of Section \ref{DM} are restricted to gradient flow systems, we will assume that both the unperturbed and perturbed systems are well approximated by stochastically forced gradient flows.

For gradient flow systems, we have $F(x)=-\nabla U(x)$ and $\delta F(x)=-\nabla \delta U(x)$, where the perturbation is a time-independent potential. Given a data set from $\{x_i\}_{i=1}^N$ sampled from $\peq$ and the functional form of $\delta U$, it is tempting to proceed as in Section~\ref{UQ} by writing $p^\delta(x,t) = \sum_i c_i(t) \psi_i(x)$ and projecting \eqref{responseeq} to eigenfunctions of $\mathcal{L}^*$. The main issue with this projection is that it is difficult to computationally attain the external perturbed forcing term, 
\[\langle\delta \mathcal{L}^*\psi_i,\psi_j\rangle_{1/\peq} = -\langle \textup{div}(\delta F\psi_i),\psi_j\rangle_{1/\peq}=\langle \delta F, \psi_i\nabla\varphi_j)\rangle = \langle \delta U, \textup{div}(\varphi_i\nabla\varphi_j)\rangle_{\peq}.\]
Even if $\delta F$ or $\delta U$ are known, we have no provably convergent technique to compute the gradient or divergence operators on the unknown Riemannian manifold. In fact, the technique discussed in Section \ref{DM} does not provide a framework to represent vector fields on manifolds.

As a remedy, we propose to project \eqref{responseeq} onto the eigenfunctions of $\tilde{\mathcal{L}}^*$ which can be accessed from the functional form of $\delta U$ and the data set $\{x_i\}$. This is an application of a result in \cite{VB}, which showed that one can approximate the generator of any gradient flow system with a (partially) known potential even if the data is sampled from a different gradient flow system. In particular, it was shown (see eq.(17) in \cite{VB}) that for general choice of bandwidth function $\rho$ in \eqref{VBkernel} and any $f\in L^2(\mathcal{M},\peq)\cap\mathcal{C}^3(\mathcal{M})$, the matrix $L_\epsilon$ constructed in \eqref{approxgenerator} converges to a weighted Laplacian in the following sense,  
\begin{align}
L_\epsilon f = \Delta f+2\left(1+\frac{d}{4}\right)\nabla f\cdot \frac{\nabla \peq}{\peq} + (d+2) \nabla f\cdot \frac{\nabla\rho}{\rho}+\mathcal{O}(\epsilon).\label{Lgeneral}
\end{align}
If we choose the following bandwidth function,
\begin{align}
\rho = \peq^{-1/2} e^{-\frac{\delta U}{D(d+2)}},\label{modbandwidth}
\end{align}
where $\peq$ is estimated from the data set as before, substituting this into \eqref{Lgeneral}, we obtain,
\begin{align}
L_\epsilon f &= \Delta f+\left(2+\frac{d}{2}\right)\nabla f\cdot \frac{\nabla \peq}{\peq}  - \nabla f\cdot \frac{\nabla\delta U}{D} -\left(\frac{d}{2}+1\right)\nabla f\cdot \frac{\nabla\peq}{\peq} + \mathcal{O}(\epsilon) \nonumber\\ &= \Delta f +\nabla f\cdot \frac{\nabla \peq}{\peq} - \nabla f\cdot \frac{\nabla\delta U}{D} + \mathcal{O}(\epsilon)\nonumber\\ &= \Delta f - \nabla f\cdot \frac{(\nabla U +\nabla\delta U)}{D} + \mathcal{O}(\epsilon)\nonumber\\ &= D^{-1}\tilde{\mathcal{L}}f+ \mathcal{O}(\epsilon).\label{newL}
\end{align}
With this result, we construct a stochastic matrix $L_\epsilon$ following exactly the procedure described in Section~\ref{findingBK}, except with the bandwidth function in \eqref{modbandwidth}. 

Let $\tilde{\lambda}_i$ and $\tilde\varphi_i$ be the eigenvalues and eigenfunctions of $L_\epsilon$, which is a discrete approximation to the continuous operator $D^{-1}\tilde{\mathcal{L}}$ as shown in \eqref{newL}. 
Following the same argument as in Section~\ref{UQ}, the corresponding Fokker-Planck operator, $\tilde{\mathcal{L}}^*$, has eigenvalues $D\tilde{\lambda}_i$ with eigenfunctions, $\tilde\psi_i = \tilde\varphi_i \tilde{p}_\textup{eq}$, where $\tilde{p}_\textup{eq} \propto\peq e^{-\delta U/D} $ is the equilibrium measure of the perturbed system. Now letting $p^\delta(x,t) = \sum_i \tilde{c}_i(t) \tilde{\psi}_i(x)$ and projecting \eqref{responseeq} with initial condition, $p^\delta(x,0)=\peq(x)$, on the perturbed coordinate basis, the evolution of the perturbed density becomes a linear forecasting problem,
\begin{align}\label{responseprojection} 
\frac{d\tilde{c}_i}{dt} = D\tilde\lambda_i \tilde c_i, \quad\quad
\tilde c_i(0) = \langle \peq, \tilde\psi_i \rangle_{1/\tilde{p}_\text{eq}}.  
\end{align}
Therefore, the response in \eqref{response} formula can be approximated by,
\begin{align}\label{responsefunctionalexpectation}
\delta\mathbb{E}[A](t) \approx \sum_{i=1}^M (\tilde{c}_i(t)-\tilde{c}_i(0)) \langle A, \tilde\varphi_i \rangle_{\tilde{p}_\text{eq}}, 
\quad t\geq0\end{align}
utilizing \eqref{functionalexpectation} on the perturbed diffusion coordinates. Notice that the zeroth mode is excluded in the summation in \eqref{responsefunctionalexpectation} since $\tilde{\lambda}_0=0$ and therefore, $\tilde{c}_0(t)=\tilde{c}_0(0)$, for any $t$.

A complication in numerically implementing the above technique is that the eigenfunctions $\tilde \varphi_i$ are orthonormal with respect to the inner product $\langle\cdot,\cdot\rangle_{\tilde p_{\rm eq}}$ weighted by the perturbed equilibrium density, $\tilde p_{\rm eq}$, but we can only estimate inner products weighted with respect to the unperturbed equilibrium density, $\peq$.  This requires us to rewrite all of the inner products with respect to $\peq$. In particular, we can rewrite the initial coefficients in \eqref{responseprojection} as,
\[ \tilde c_i(0) = \langle \peq,\tilde \psi_i \rangle_{1/\tilde p_{\rm eq}} = \left<\peq,\tilde \varphi_i \right> = \langle 1,\tilde\varphi_i \rangle_{\peq} \approx \frac{1}{N} \sum_{j=1}^N \tilde\varphi_i(x_j).\]
With these initial conditions, we have explicit solutions for the linear problem in \eqref{responseprojection} given by $\tilde c_i(t)= e^{D\tilde\lambda_it}\tilde c_i(0)$.  To compute the response formula in \eqref{responsefunctionalexpectation}, we evaluate the inner product terms in \eqref{responsefunctionalexpectation} as follows,
\[\langle A, \tilde\varphi_i\rangle_{\tilde{p}_\text{eq}} = \langle A,  \frac{\tilde{p}_\text{eq}}{\peq} \tilde\varphi_i \rangle_{\peq} \approx \frac{1}{N}\sum_{j=1}^N A(x_j) \frac{\tilde{p}_\text{eq}(x_j)}{\peq(x_j)} \tilde\varphi_i(x_j) = \frac{1}{NZ}\sum_{j=1}^N A(x_j) e^{-\delta U(x_j)/D}\tilde\varphi_i(x_j),\]
using the fact that the equilibrium density of the perturbed system is given by, $\tilde p_{\rm eq} = \frac{1}{Z} \peq e^{-\delta U/D}$, where $Z$ is the normalization factor for the density $\tilde p_{\rm eq}$ obtained as discussed in Section~\ref{QoI}. 

\section{Numerical Results}\label{sec5}

For a state variable $x$ evolving according to \eqref{SDE}, we have shown how to estimate the density function $p(x,t)$, which solves either the forecasting, filtering, or response problem.  We have also shown how to find the expectation of any quantity of interest, $A(x)$, with respect to this measure.  In this section, we first validate the numerical algorithms developed in Sections~\ref{UQ} and \ref{sec4} on a simple linear gradient flow system.  We then test our algorithms on a chaotically forced gradient flow system which is known to be well approximated as a stochastically forced gradient flow. In each example we will use only $M=50$ eigenvectors in the non-parametric model.

\subsection{Linear Example: Ornstein-Uhlenbeck processes}\label{sec51}
We first validate our algorithms on the Ornstein-Uhlenbeck processes on $\mathcal{M}=\mathbb{R}$, given by,
\begin{align}\label{OU} dx = -x\, dt + \sqrt{2}\, dW_t, \end{align}
where $dW_t$ denotes standard Gaussian white noise.  Our technique assumes that we do not know the model in \eqref{OU}, but instead we are only given sample solutions of \eqref{OU} as a training data set.  In our numerical experiment, we generate a training data set consisting of $N=10000$ data points $x_l = x(l\Delta t)$ by numerically integrating \eqref{OU} with the Euler-Maruyama method and sub-sample at every $\Delta t=0.2$ time units.  We intentionally sub-sample with a large $\Delta t$ so that the data set is approximately independently sampled from the invariant measure. We numerically validate our methods on this simple system \eqref{OU} since it is possible to find an analytic expression for the evolution of the moments (see Appendix \ref{OUanalysis}) for the forecasting and response problems. For the filtering problem, we compare the results to the Kalman-Bucy solutions which are optimal in a least square sense \cite{kalman:61}.  

\subsubsection{Forecasting}\label{forecastalgorithm}

\begin{figure}[h]
\centering
\includegraphics[width=0.45\textwidth]{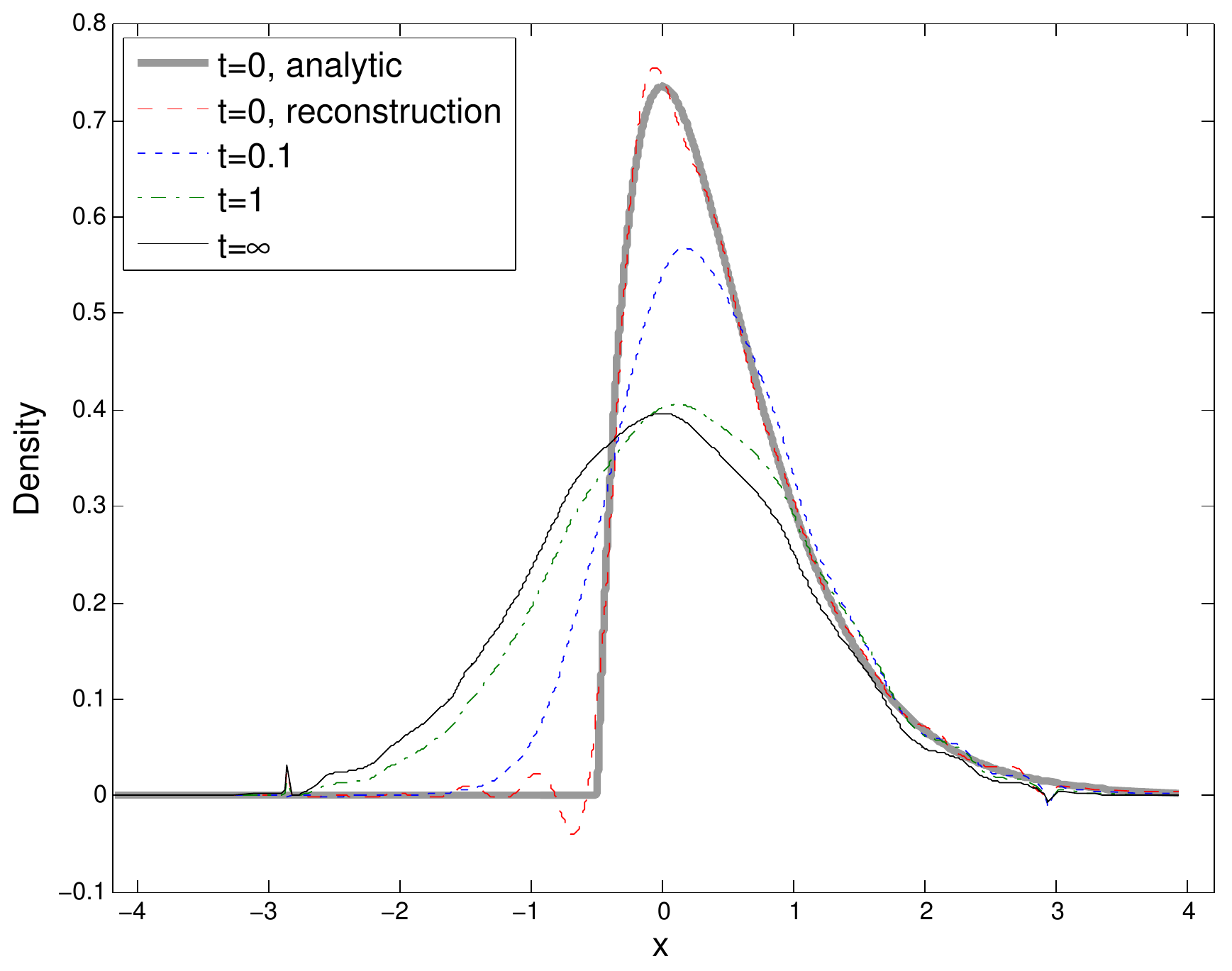}\includegraphics[width=0.45\textwidth]{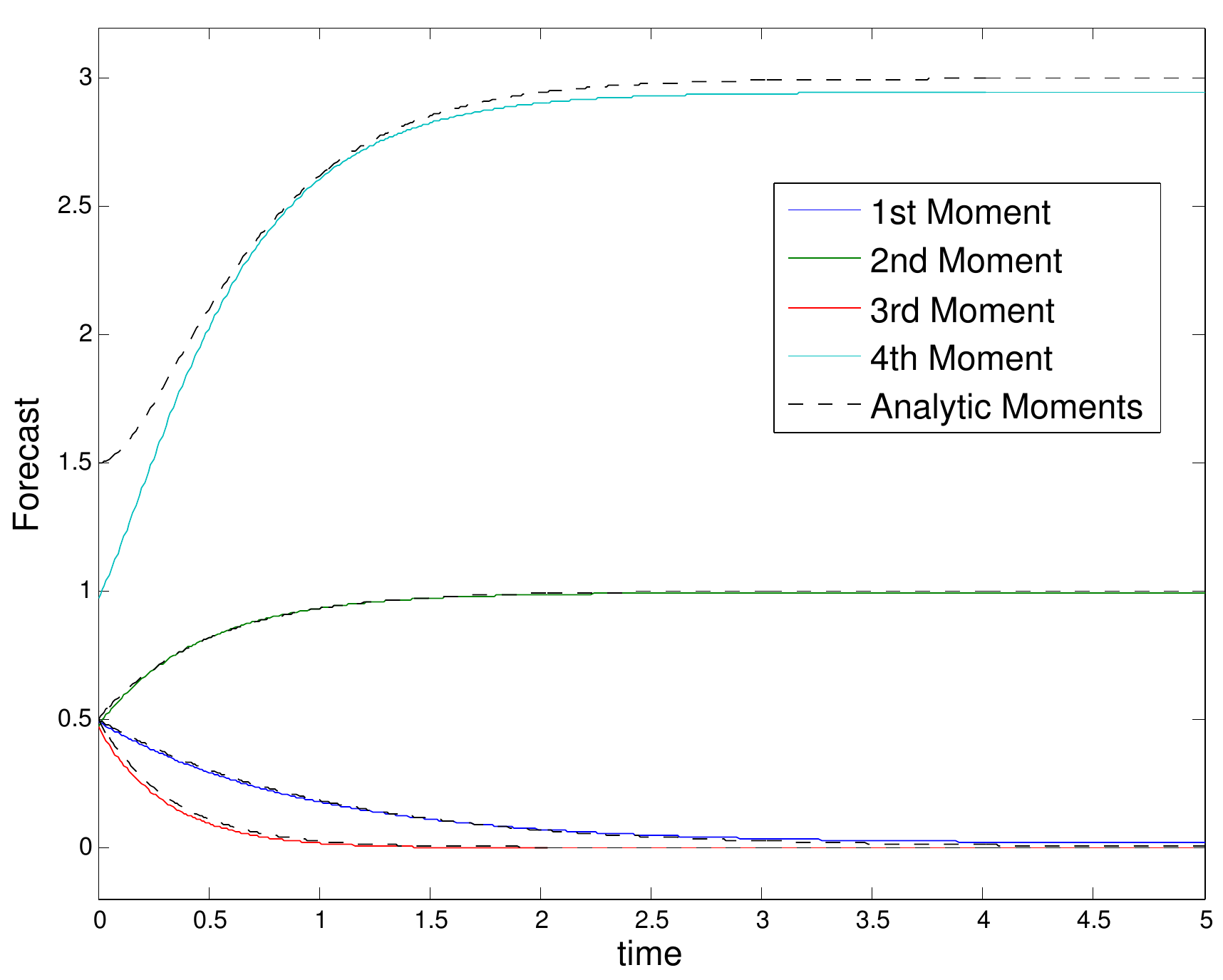}
\caption{\label{forecastingfig} Forecasting the evolution of a Gamma distributed initial condition in the Ornstein-Uhlenbeck system.  Left: Evolution of the density. Right: Evolution of the moments.}
\end{figure}

To validate our forecasting algorithm, we take the initial condition to be a Gamma distribution with density $p_0(x) = 4 (x+1/2) e^{-2(x+1/2)}$ on $[-1/2,\infty)$. For this distribution the first four moments are, $\bar x = \mathbb{E}_{p_0}[x] = 1/2$, $\mathbb{E}_{p_0}[(x-\bar x)^2] = 1/2$, $\mathbb{E}_{p_0}[(x-\bar x)^3] = 1/2$, and $\mathbb{E}_{p_0}[(x-\bar x)^4] = 3/2$.  We assume that the initial distribution is known, and we simply evaluate $p_0$ on the training data set to find $(\vec p_0)_l = p_0(x_l)$.  The results of the forecasting algorithm are shown in Figure \ref{forecastingfig} where we reconstruct the forecast density at various times as it approaches the invariant measure. Notice the oscillation in the reconstructed $p_0$ near the non-smooth part is a Gibbs phenomenon-like behavior. We also show the evolution of the first four centered moments which are evaluated as linear combinations of the uncentered moments $A(x) = x^r$, $r=1,\ldots, 4$. Notice that the initial kurtosis of the Gamma distribution is difficult to estimate since the data on the tail is very sparse (there are only 3 sample points for $x>3$ in the training data set).  The Monte-Carlo integral converges slowly since the Gamma distribution decays slower than the Gaussian distribution.

\subsubsection{Filtering}\label{filteralgorithm}

To verify our filtering algorithm, we consider the observation \eqref{obseq} with $h(x) = x$ so that,
\[ dz = x\, dt + \sqrt{R}\,dW_t. \]
The filtering problem for \eqref{OU} with this observation has optimal solutions given by the Kalman-Bucy equations \cite{kalman:61}, which we approximate for a finite observation time step $\Delta t$ by the discrete time Kalman filter.  Since discrete time observations are typically given in the form,
\[ \mathcal{Z}(t_i) = x(t_i) + \sqrt{R_o}\,\omega(t_i) \]
where $\omega(t_i)$ are independent random samples from $\mathcal{N}(0,1)$, we will approximate $dz$ for a finite time step $\Delta t = 0.01$ as,
\[ dz(t_i) \approx x(t_i)\Delta t + \sqrt{\Delta t R_o}\, dW_{t_i} = \mathcal{Z}(t_i)  \Delta t . \]
So in the numerical algorithm of Section \ref{nonlinearfilter} we use $dz(t_i) = \mathcal{Z}(t_i)\Delta t$ and we set $R = \Delta t R_o$.  

\begin{figure}[h]
\centering
\includegraphics[width=0.45\textwidth]{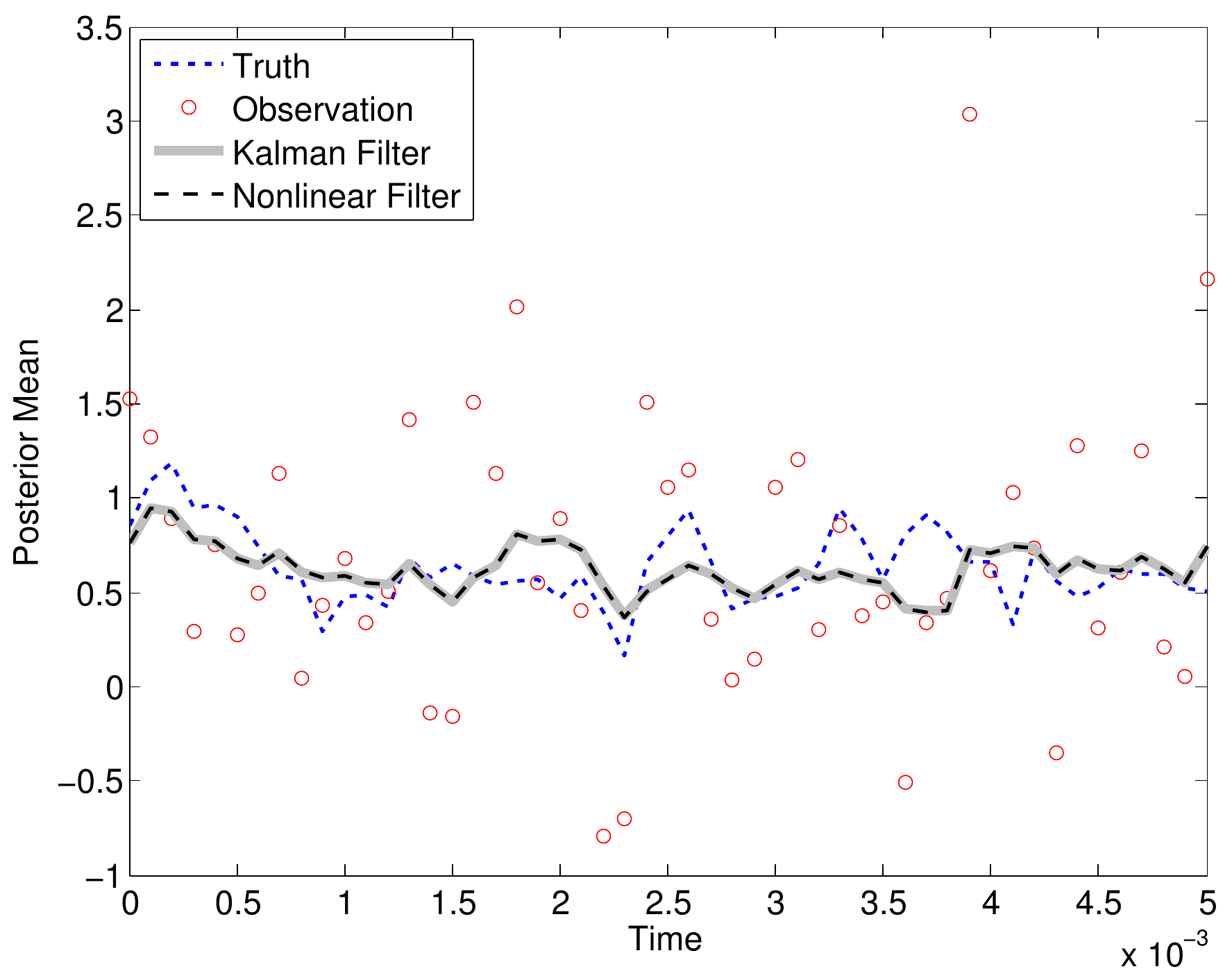} \includegraphics[width=0.45\textwidth]{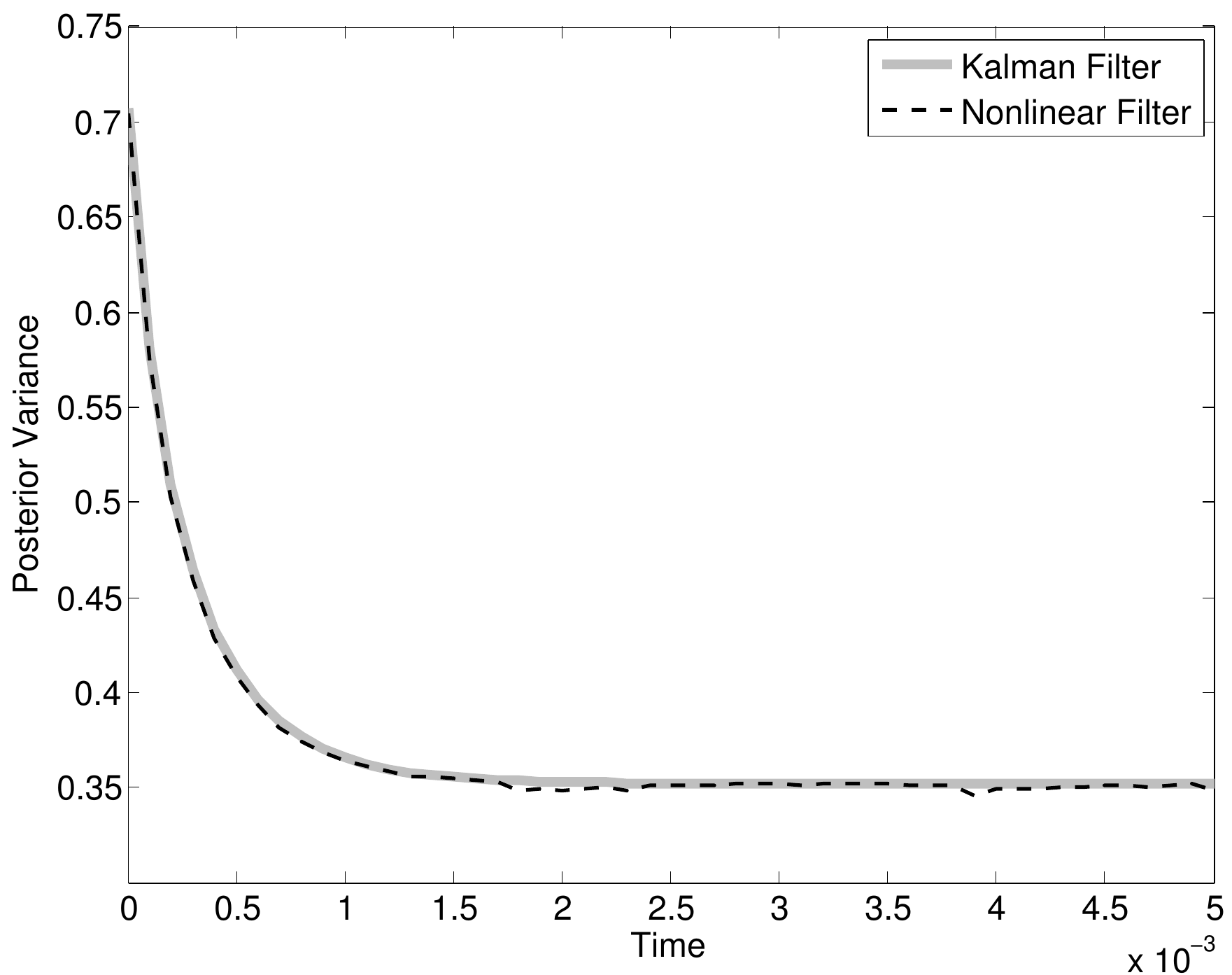}
\caption{\label{filterfig} Nonlinear filter mean estimates as functions of time, compared to the truth, Kalman filter solutions, and observations (left). The posterior covariance estimates as functions of time (right). The nonlinear filter recovers the Kalman filter solutions in this linear example.  The posterior evolution is shown in {\tt fig2.mov}.}
\end{figure}

In our numerical simulation, we take $R_o=1$ and assume an initial distribution $p_0(x)=\peq(x)$, which is the Gaussian invariant measure of \eqref{OU}.  The results are shown in Figure \ref{filterfig}, where we show that we recover the same statistical solutions (mean and covariance) as the discrete Kalman filter.  We also include a video, {\tt fig2.mov}, comparing the evolution of the reconstructed posterior density with the Kalman filter posterior.

\subsubsection{Response}\label{responsealgorithm}

\begin{figure}[h]
\centering
\includegraphics[width=0.5\textwidth]{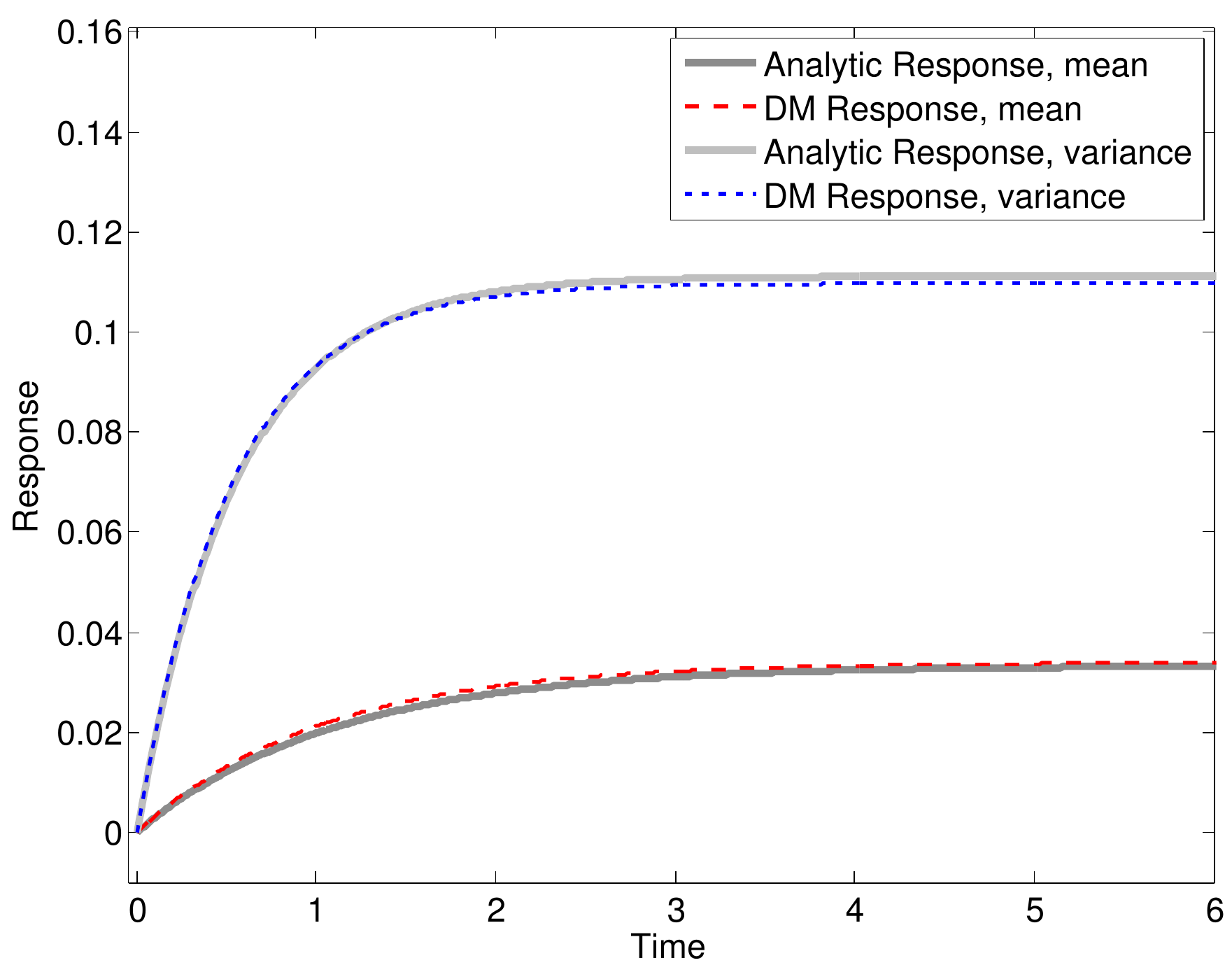}
\caption{\label{responsefig} Nonlinear response of the first two moments of the Ornstein-Uhlenbeck system to a perturbation of the potential function.  The analytic response is derived in Appendix \ref{OUanalysis} and the DM response is our nonparametric estimate.}
\end{figure}

To validate the response algorithm of Section \ref{sec4}, we perturb the Ornstein-Uhlenbeck system in \eqref{OU}, which has potential $U(x) = x^2/2$, with a potential $\delta U(x) = \frac{a}{2}(x-b/a)^2$.  The vector field for the potential $U+\delta U$ becomes $-\nabla(U+\delta U) = -x-a x + b$ so that the perturbed system is again an Ornstein-Uhlenbeck system with different mean and damping parameters.  This fact allows us to easily compute the analytic response (see Appendix \ref{OUanalysis}) which we compare to the numerical estimate using the above algorithm in Figure \ref{responsefig}. In this example we choose $a=-0.1$ and $b=0.03$, since the damping is decreased by this perturbation the variance increases in response while the effect of $b$ is a shift in the mean. Notice the accuracy of the response estimates for the mean and variance at all times.

\subsection{Nonlinear Example: Chaotically driven double well potential}\label{sec52}

In this section we test our nonparametric UQ methods on the following nonlinear system,
\begin{align}
\label{doublewellfull}
\begin{split}
\dot x &= x-x^3 + \frac{\gamma}{\epsilon} y_2,\\
\dot y_1 &=\frac{10}{\epsilon^{2}}(y_2 - y_1), \\
\dot y_2 &=\frac{1}{\epsilon^{2}}(28 y_1 - y_2 - y_1 y_3),\\
\dot y_3 &=\frac{1}{\epsilon^{2}}(y_1 y_2 - \frac{8}{3}y_3).
\end{split}
\end{align}
which was discussed in \cite{gks:04} as an example of a stochastic homogenization problem and used in \cite{mg:12} as a test model for a reduced filtering problem. The fast variables $(y_1,y_2,y_3)$ solve a chaotic Lorenz-63 system, and $x$ is a scalar ODE driven by this chaotic oscillator with characteristic time $\epsilon^2$. In \cite{gks:04,mg:12} it was shown that for $\gamma, \epsilon$ sufficiently small the dynamics of the variable $x$ are well approximated by the reduced stochastic model.
\begin{align}\label{doublewellreduced} 
dX = X(1-X^2) \, dt + \sigma \, dW_t,
\end{align}
where $\sigma$ is a diffusion constant that can be numerically estimated by a linear regression based algorithm or by matching the correlation time of $\gamma y_2$.  This reduced model is a stochastically forced gradient flow system, where the potential function $U(X) = X^2/2 - X^4/4$ is a standard double well potential.  While this suggests that our technique should be successful, we emphasize that our method makes no use of either \eqref{doublewellfull} or \eqref{doublewellreduced}, instead it relies on a training data set to represent the reduced model for the variable $x$, with no assumed parametric form.  Following \cite{mg:12}, we set $\gamma = 4/90$ and $\epsilon=\sqrt{0.1}$.  We note that \cite{mg:12} suggested $\sigma^2 =0.113$ for the homogenization and $\sigma^2 = 0.126$ for the filtering problem.

\begin{figure}[h]
\centering
\includegraphics[width=0.45\textwidth]{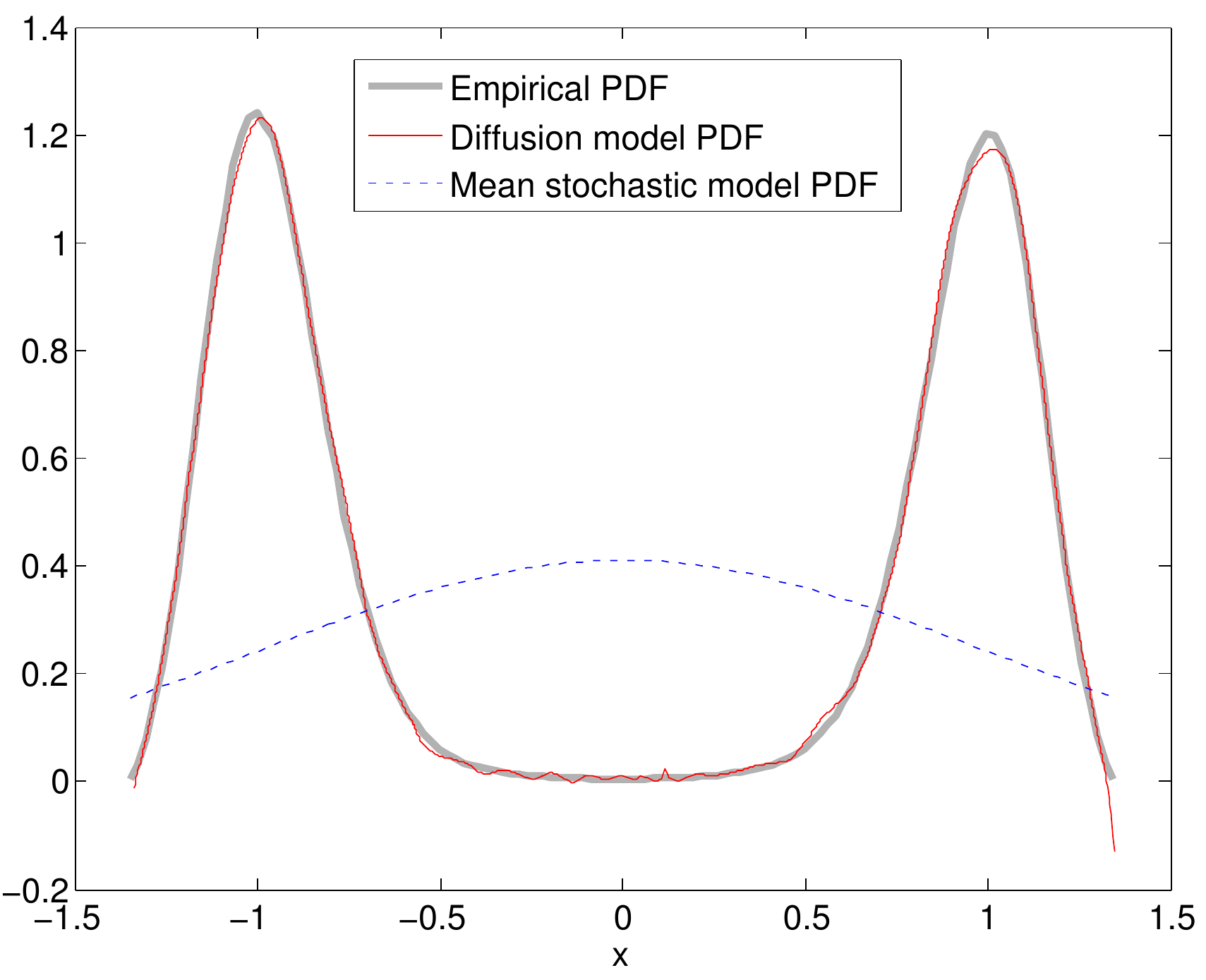}\includegraphics[width=0.45\textwidth]{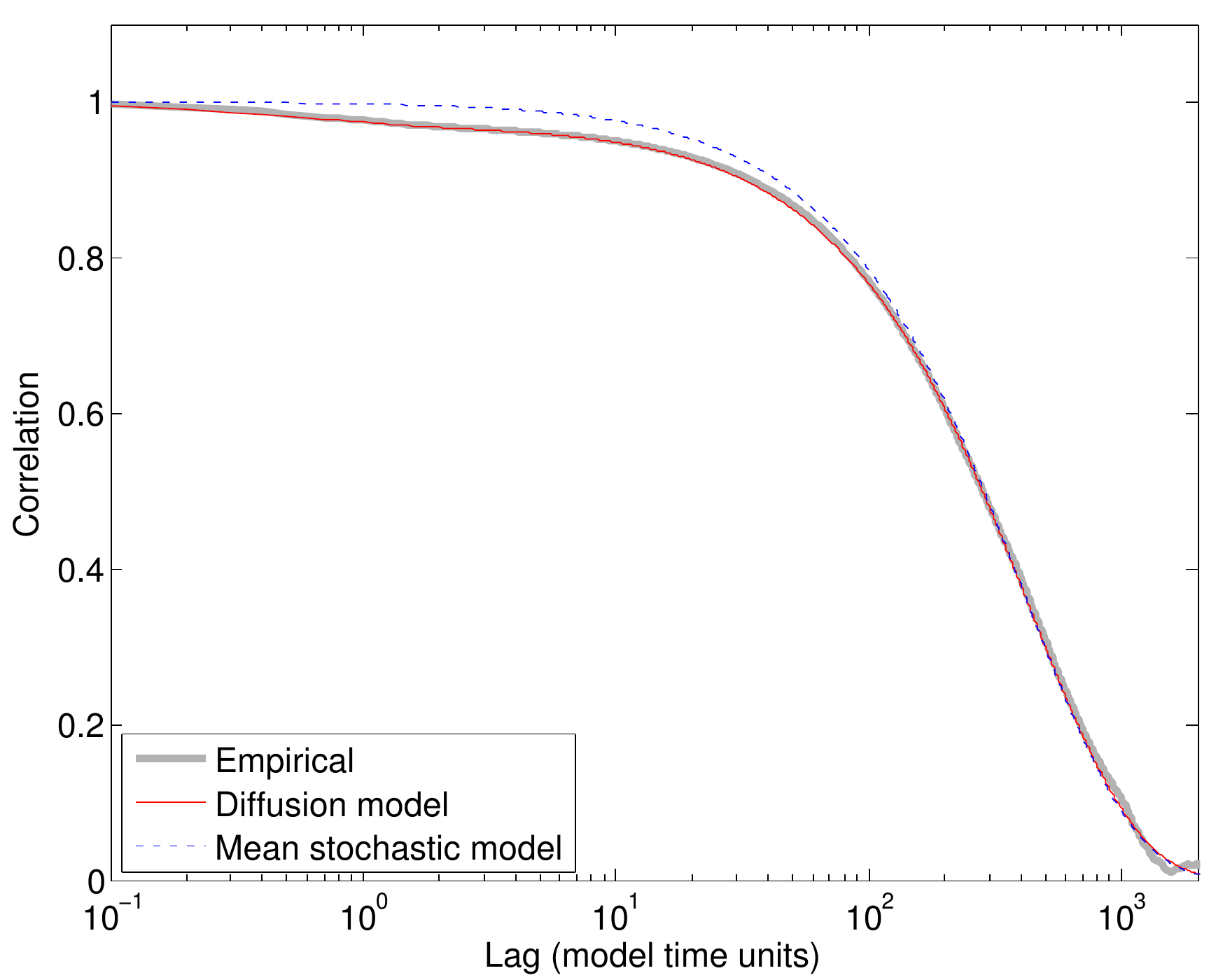}
\caption{\label{statistics} Left: Invariant measure estimated by the kernel. Right: Tuning the stochastic forcing using the correlation time.}
\end{figure}

In order to apply our method to this problem, we generated a training data set containing 200000 data points $\{x_j = x(t_j)\}$ with time spacing $\Delta t = t_{j+1}-t_j = .1$ by solving \eqref{doublewellfull} using a fourth order Runge-Kutta scheme with integration time step $\delta t = 0.002$.  We applied the variable bandwidth diffusion kernel as described in Section \ref{findingBK} to the data set $\{x_{10j}\}_{j=1}^{20000}$.  We then estimate the diffusion constant $D=\sigma^2/2 \approx 0.0535$ as in Section \ref{findingD}, using the entire data set $\{x_j\}$ to find the correlation time and the empirical correlation function.  In Figure \ref{statistics} we compare the invariant measure and correlation function of the nonparametric diffusion model to the empirical estimates of these statistics taken from training data set. We also show the statistics of the mean stochastic model (MSM) of \cite{mgy:10,mh:12}, which is a linear (and Gaussian) model constructed by matching the correlation time and the equilibrium density variance. From Figure~\ref{statistics}, one can see that the density and the autocorrelation function from the diffusion map fitting are very accurate.

\subsubsection{Forecasting}\label{doublewellforecast}

We first test our forecasting algorithm for an initial condition given by a normal distribution centered at $0.5$ with variance $0.01$. Unlike the Ornstein-Uhlenbeck system in Section \ref{sec51}, we do not have an analytical expressions for the evolution the moments.  For diagnostic purposes, we will compare our forecast to a Monte-Carlo simulation. We initialize an ensemble of $100000$ points $(x,y_1,y_2,y_3)^\top \in \mathbb{R}^4$ with $x$ chosen randomly according to the initial density and triples $(y_1,y_2,y_3)^\top$ randomly selected from the training data set (so that they lie on the attractor of the chaotic system).  This ensemble is then evolved according to the system \eqref{doublewellfull} with the same Runge-Kutta integrator used to generate the training data.  For any time $t$ we generate the Monte-Carlo density using a histogram of the ensemble at time $t$, and we call the moments of the ensemble the Monte-Carlo moments.  These Monte-Carlo estimates make use of the full true model and are computationally more intensive than our model free method, and we will consider them the `truth' for the purposes of comparison.  

\begin{figure}[h]
\centering
\includegraphics[width=0.425\textwidth]{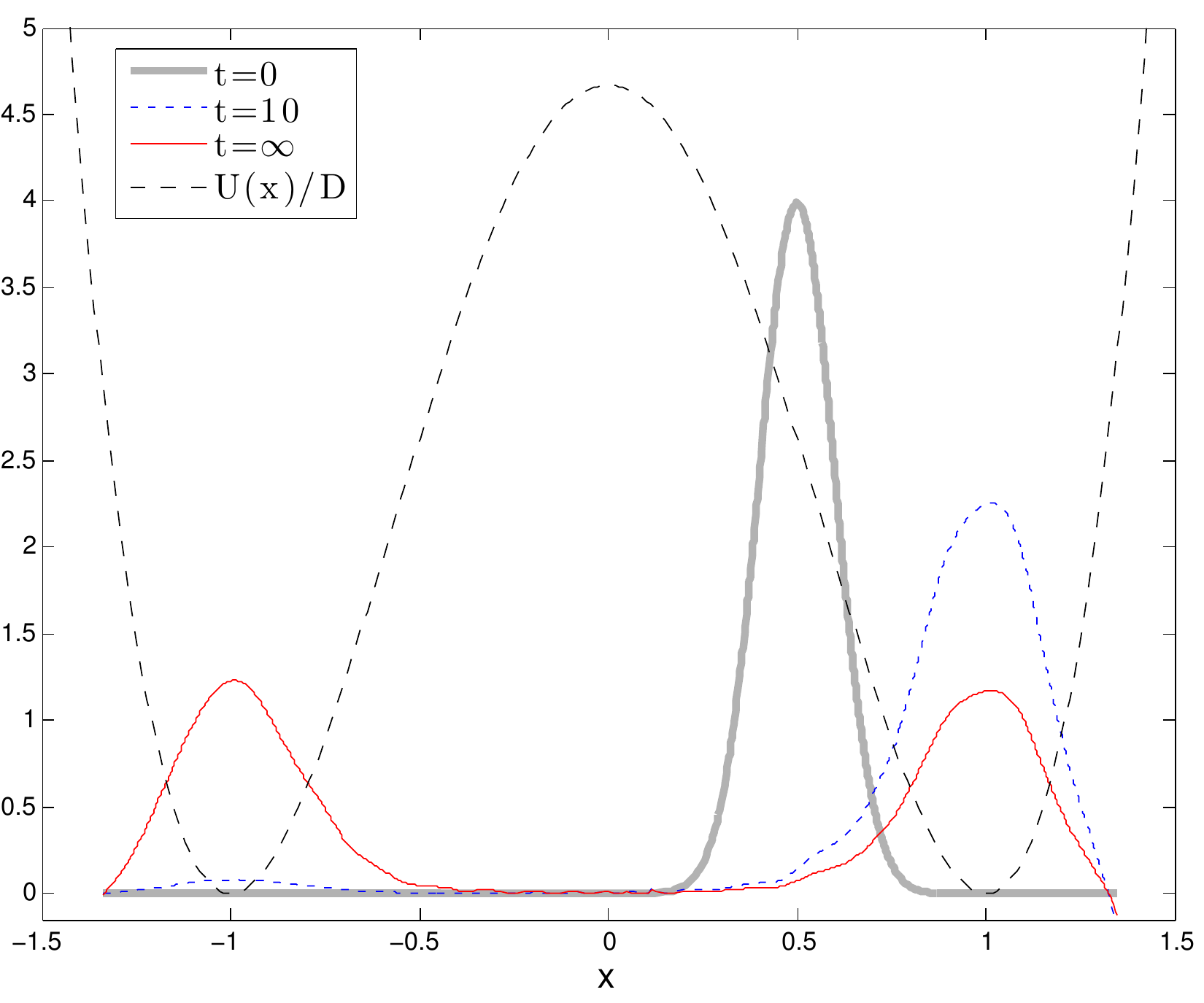}\includegraphics[width=0.45\textwidth]{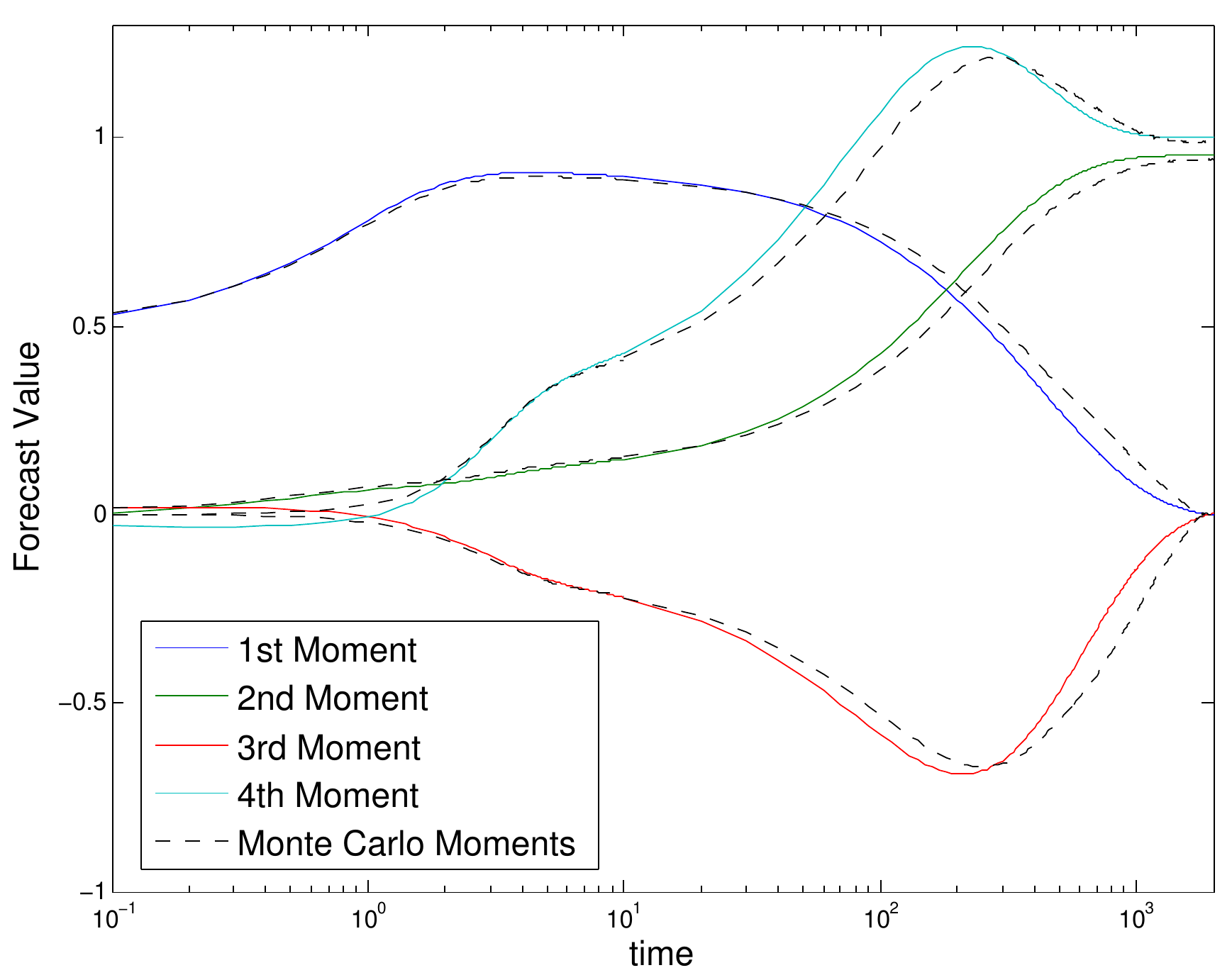}
\caption{\label{doublewellprediction} Left: Evolution of the density is summarized in snapshots for $t=0,10,\infty$ and the potential function is shown.  Right: Non-parametric estimates for the evolution of the moments is compared to the Monte-Carlo simulation.  The reconstructed density evolution is shown in {\tt fig5a.mov} and {\tt fig5b.mov}. }
\end{figure}

Next, we estimate the evolution of the density and the first four centered moments using the technique of Section \ref{UQ}.  For this initial condition, the density moves to the right and centers on the right well of the potential on the first 10 model time units. Over the next 1000 model time units, small amounts of density slowly migrate across the barrier between the two potential wells until the system is at equilibrium.  In Figure \ref{doublewellprediction} we compare the moment estimates of the nonparametric method to the Monte-Carlo moments.  We also show snapshots of the evolution of the density as reconstructed from the diffusion coordinates at times $t = 0, 10, \infty$. The evolution of the reconstructed density is compared to the Monte-Carlo simulation in {\tt fig5a.mov} up to time $t=10$ and {\tt fig5b.mov} up to time $t=2000$.

\subsubsection{Filtering}

\begin{figure}[h]
\centering
\includegraphics[width=0.45\textwidth]{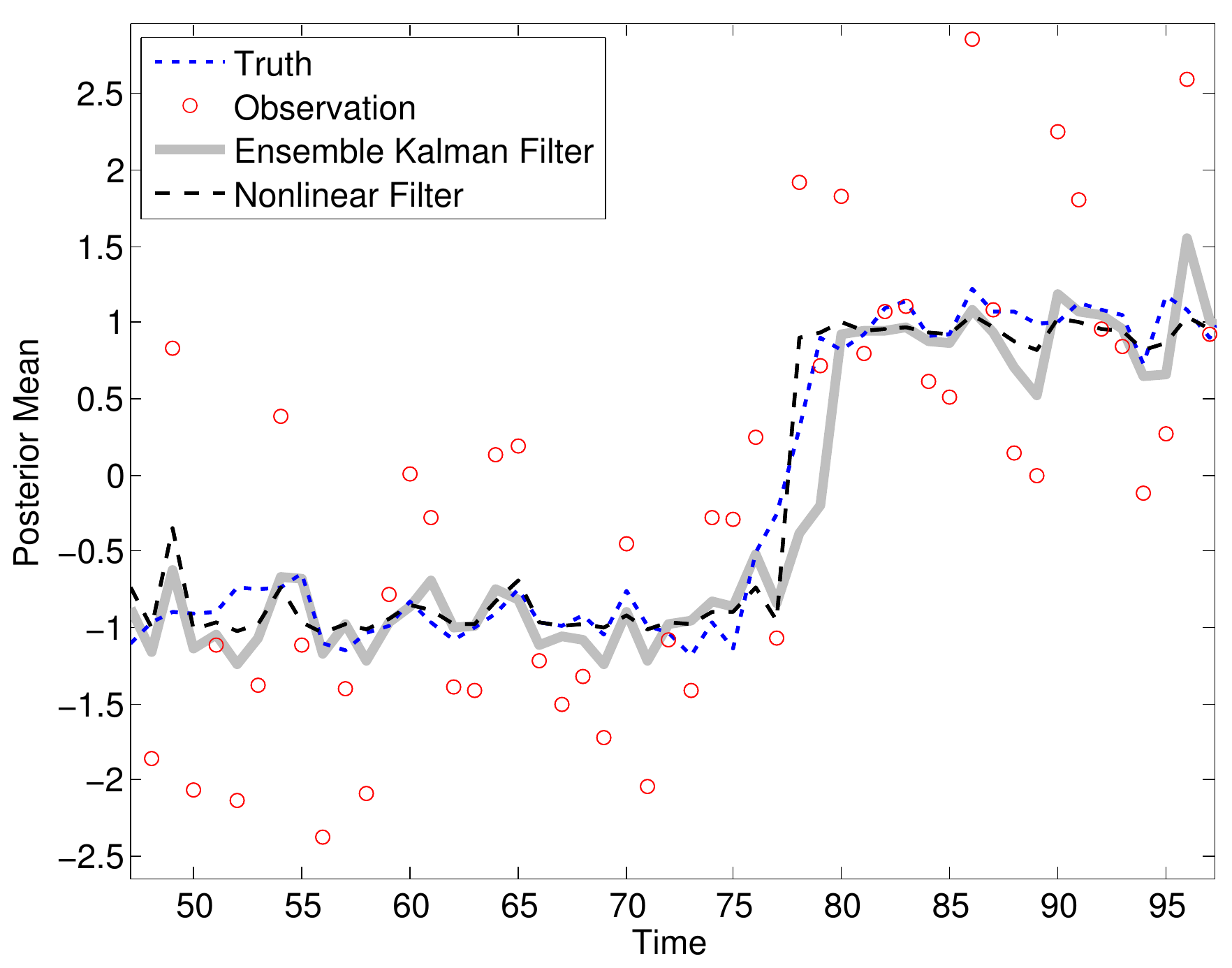}\includegraphics[width=0.45\textwidth]{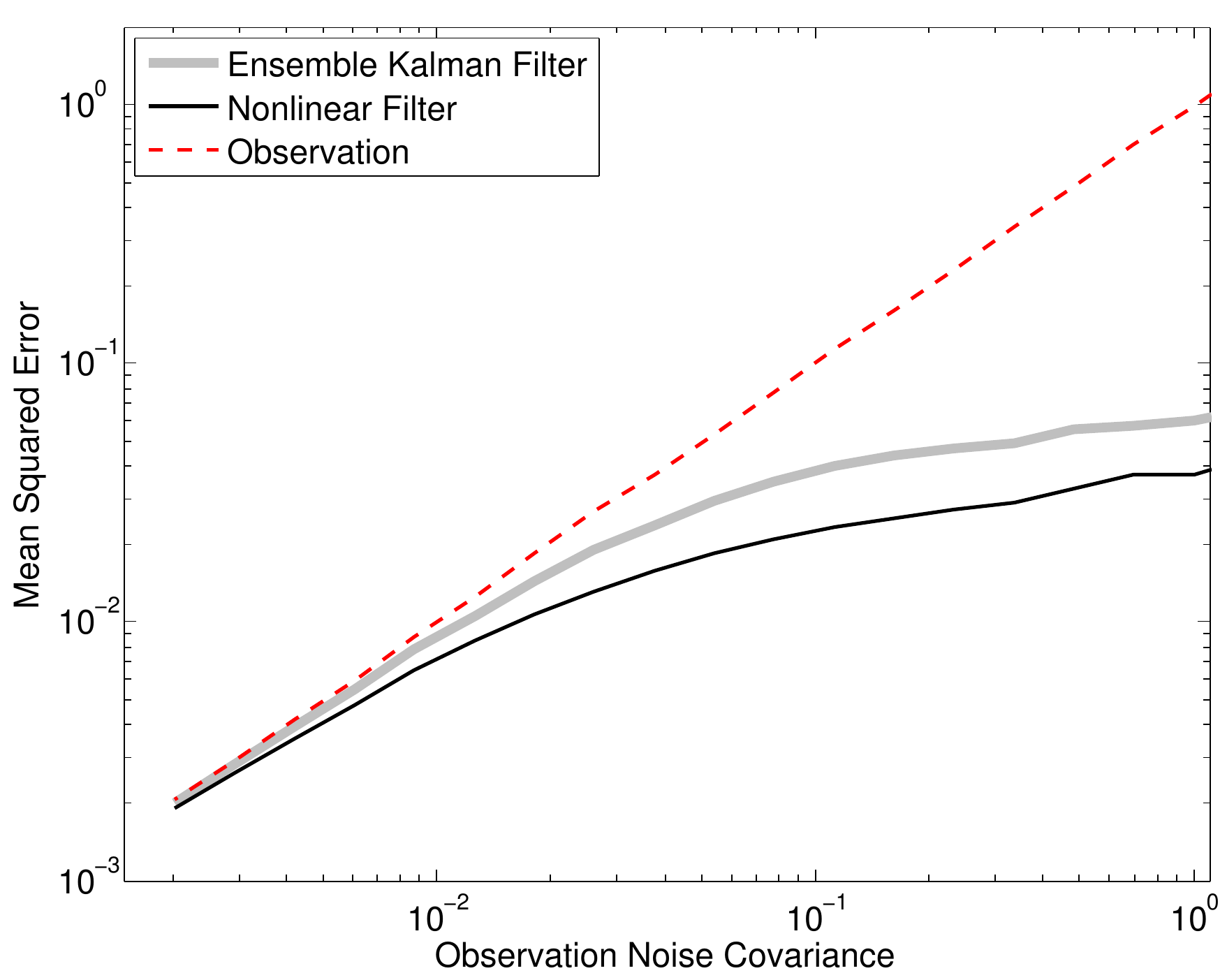}
\caption{\label{doublewellfilter} Left: Evolution of the posterior mean during a transition from the left potential well to the right potential well.  Right: Mean squared error between filter estimates and the true state as a function of the observation noise covariance. The evolution of the posterior densities is shown in {\tt fig6.mov}.}
\end{figure}

We now test our nonlinear filtering algorithm on the chaotically driven double well potential using the observation \eqref{obseq} with $h(x) = x$ for several values of observation noise variance, $R$, ranging from $10^{-3}$ to $1$.  We choose a time discretization with $\Delta t = 1$ between observations, which is relatively long compared to $\Delta t=0.1$ considered in \cite{mg:12}. For each noise level $R$, we generate 10000 discrete observations using the full system \eqref{doublewellfull}, integrated as in Section \ref{doublewellforecast} and observed with the discretization of \eqref{obseq}. We then apply the filtering algorithm from Section \ref{nonlinearfilter}, which we refer to as the nonlinear filter in Figure \ref{doublewellfilter} where we show the mean estimate as the quantity of interest. For comparison we also use the standard ensemble Kalman filter (EnKF) for the reduced model \eqref{doublewellreduced} using $\sigma^2 = .126$ as suggested in \cite{mg:12,gks:04} with ensemble size 50 without additional variance inflation. 

In Figure \ref{doublewellfilter} we see that the nonlinear filter tracks the transition between the potential wells more accurately compared to EnKF for large observation noise, $R=0.5$ (see left panel). In the large noise regime, the posterior estimate is very far from Gaussian. Conversely, for very small observation noise and also for very short discretization time, the performance of nonlinear filter and the EnKF are similar, their mean-squared errors are close to the observation noise variance (see the right panel of Figure~\ref{doublewellfilter}). In the small noise regime, the posterior density is close to Gaussian since it is very close to the observation likelihood. The evolution of the reconstructed non-Gaussian posterior density is compared to the Gaussian posterior density of the EnKF in {\tt fig6.mov}.

\begin{figure}[h]
\centering
\includegraphics[width=0.45\textwidth]{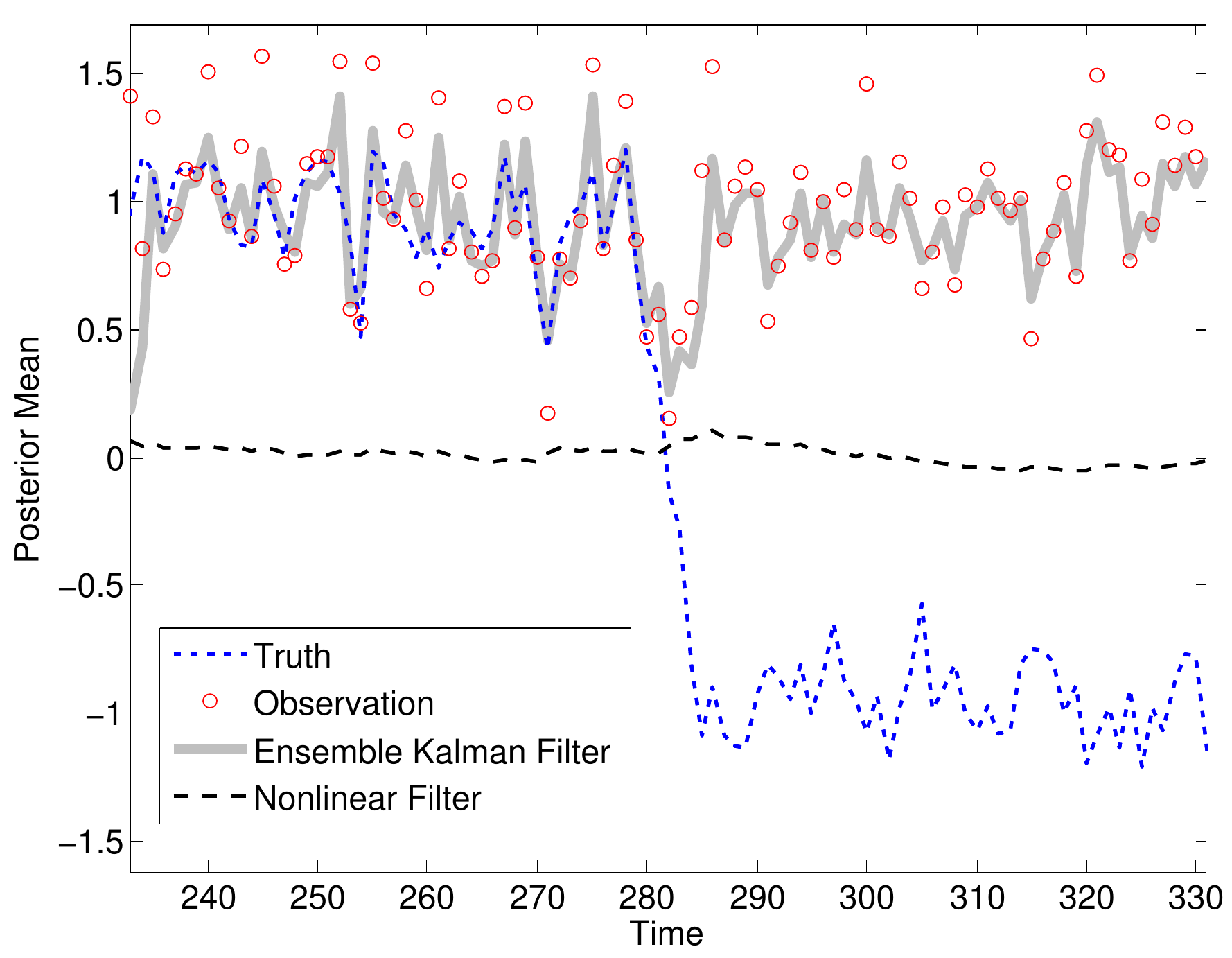}\includegraphics[width=0.45\textwidth]{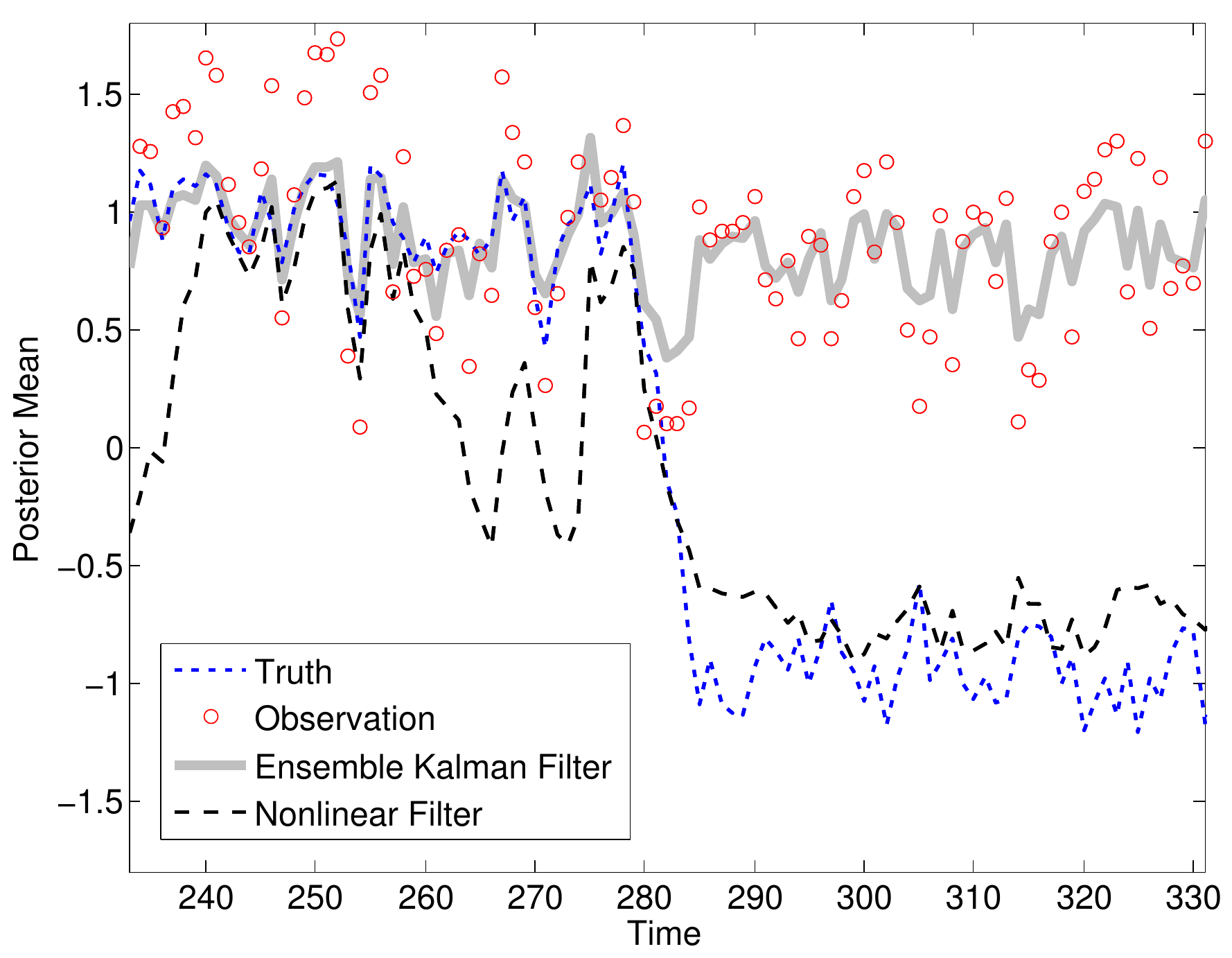}
\caption{\label{doublewellnonlinearobs} Left: Evolution of the posterior mean for $h(x) = |x|$.  Right: Evolution of the posterior mean for $h(x) = (x-0.05)^2$. Both filtering experiments assimilate observations at every $\Delta t=1$ model time unit and noise variance $R=0.05$.  The evolution of the posterior densities are shown in {\tt fig7a.mov} and {\tt fig7b.mov} respectively.}
\end{figure}

Next, we consider a nonlinear observation with $h(x) = |x|$, which is a pathological observation since it is impossible to determine from this observation which potential well the truth lies in. As shown in Figure \ref{doublewellnonlinearobs}, the nonlinear filter respects this symmetry with the posterior mean always very close to zero. In contrast, the EnKF outperforms the nonlinear filter when it happens to be on the right side of the well, however it frequently does not track the transitions leading to a higher overall mean squared error.  The evolution of the posterior densities are shown in {\tt fig7a.mov} which shows the posterior density of the nonlinear filter to be a symmetric bi-modal distribution.  We also consider another nonlinear observation with $h(x) = (x-0.05)^2$, which is nearly pathological, however, since the symmetry is not perfect, repeated observations allow the filter to determine which well the true state is in.  As shown in Figure \ref{doublewellnonlinearobs} and in {\tt fig7b.mov}, the posterior of the nonlinear filter becomes increasingly bi-modal when observations occur near $0.05$.  However, with repeated observations which are far from $0.05$ the posterior density accumulates information on the true location and the incorrect mode is progressively damped away.

\subsubsection{Response}\label{dwresponsesection}

\begin{figure}[h]
\centering
\includegraphics[width=0.425\textwidth]{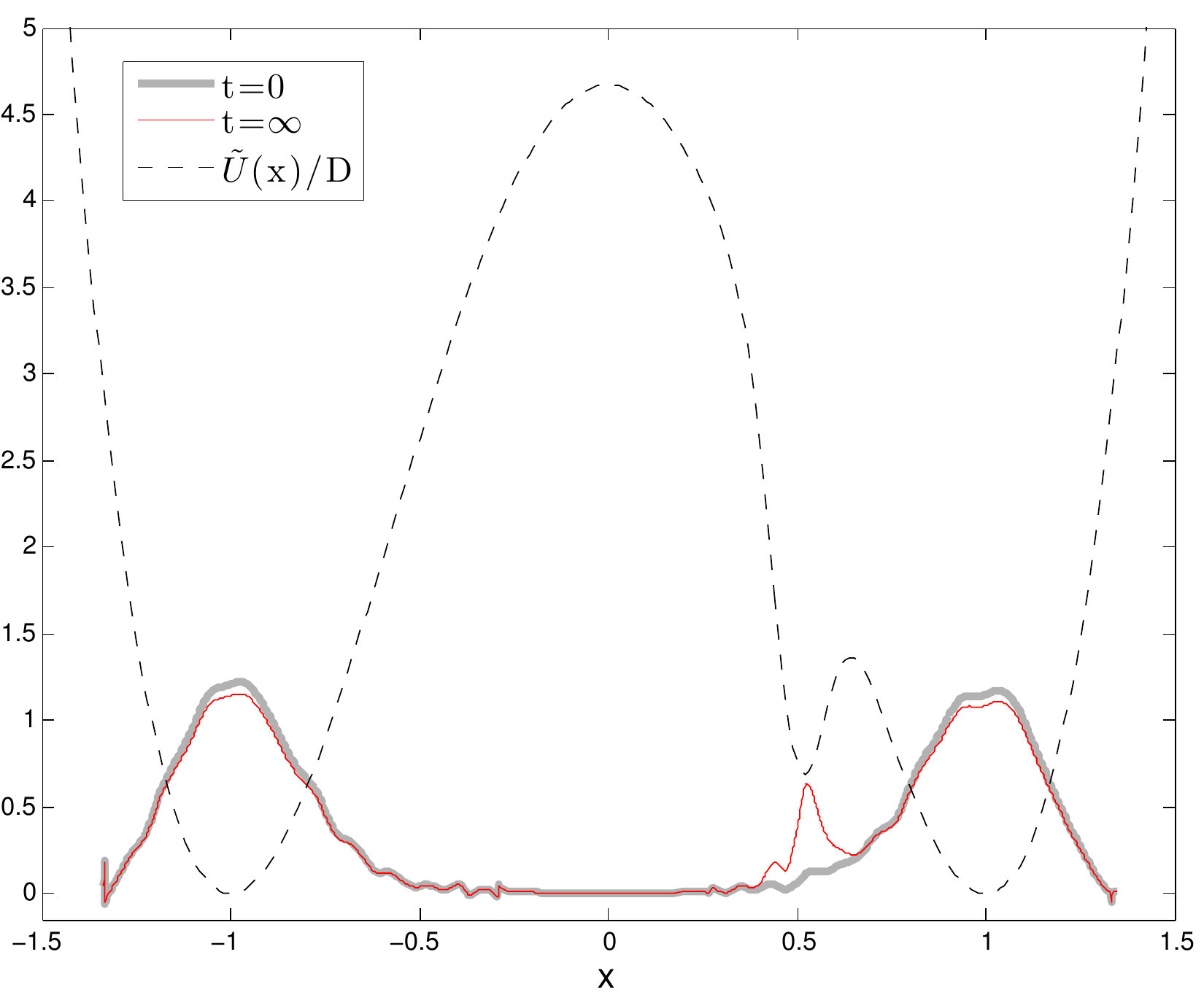}
\includegraphics[width=0.45\textwidth]{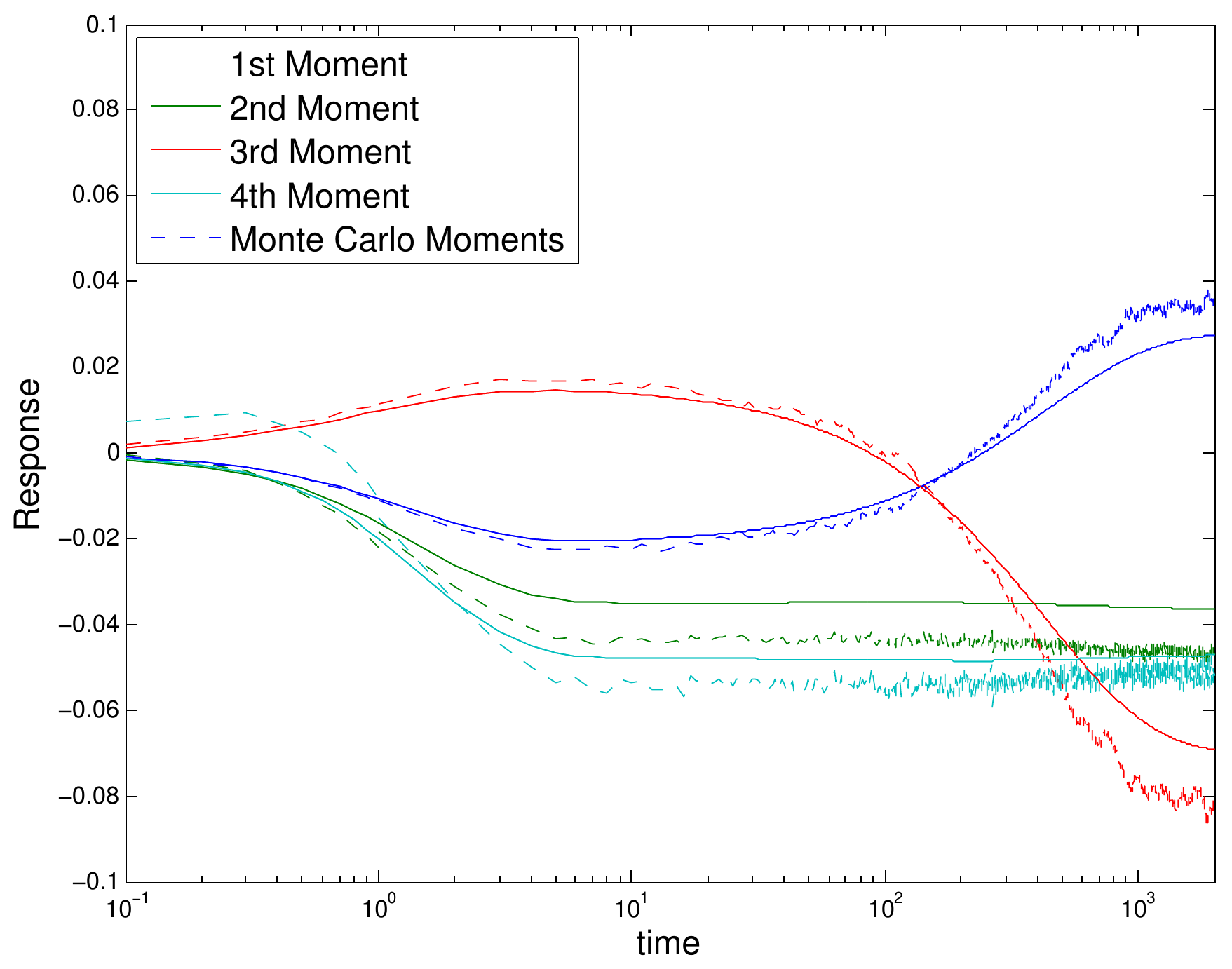}
\caption{\label{doublewellresponse} Left: Evolution of the density and the perturbed potential function.  Right: Evolution of the moments compared to those of the stochastic Monte-Carlo simulation of the reduced model with $\sigma^2=0.113$.}
\end{figure}

In this section we apply the method of Section \ref{sec4} to quantify the response of the system \eqref{doublewellfull} to perturbing the potential $U(x) = x^2/2-x^4/4$ with 
\[ \delta U(x) = -\exp(-100(x-0.5)^2)/10.\]  
This introduces a third potential well centered at $x=0.5$ as shown in Figure \ref{doublewellresponse} and changes the vector field for the variable $x$ in \eqref{doublewellfull} which becomes,
\begin{align}\label{perturbeddoublewell}
\dot x &= x-x^3 - 20(x-0.5)\exp(-100(x-0.5)^2) + (\gamma/\epsilon) y_2. 
\end{align}
We first create a benchmark for comparison by estimating the response through a Monte-Carlo simulation.  Since the system is assumed to be initially at equilibrium for the unperturbed system \eqref{doublewellfull}, we initialize an ensemble of 100000 members by simply taking the final ensemble from the long simulation in Section \ref{doublewellforecast}.  We first use a stochastic Monte-Carlo simulation by replacing $(\gamma/\epsilon)y_2$ in \eqref{perturbeddoublewell} with the stochastic term $\sigma dW_t$ with $\sigma^2 = 0.113$ as suggest in \cite{mg:12}.  This yields a reduced model as suggested in \cite{mg:12},
\begin{align}\label{reducedperturbeddoublewell}
dx &= \left(x-x^3 - 20(x-0.5)\exp(-100(x-0.5)^2)\right)\, dt + \sigma \, dW_t 
\end{align}
We apply the Euler-Maruyama scheme with integration time step $\delta t = 0.002$ to evolve this ensemble using \eqref{reducedperturbeddoublewell}.  At any fixed time $t$ we can then use the centered moments of the ensemble to produce a Monte-Carlo estimate of the response.  We also tested the reduced model using our estimate of $D$ to define $\sigma^2 = 2D= 0.107$  and the results were similar the those shown below.  

We then apply the algorithm in Section \ref{sec4} to estimate the response.  This technique uses the same training data as Section \ref{doublewellforecast} along with the perturbation $\delta U$ given above to estimate the response without any parametric model.  In the first 10 time units, density migrates quickly from the potential well centered at $x=1$ into the new potential well centered at $x=0.5$.  Over the next 1000 time units, density migrates from the potential well centered at $x=-1$ into the wells at $x=0.5$ and $x=1$, which are now more stable due to the perturbation.  In Figure \ref{doublewellresponse} we summarize the evolution of the density and compare the response estimates of the nonparametric technique to the Monte-Carlo estimates of \eqref{reducedperturbeddoublewell}.  

\begin{figure}[h]
\centering
\includegraphics[width=0.45\textwidth]{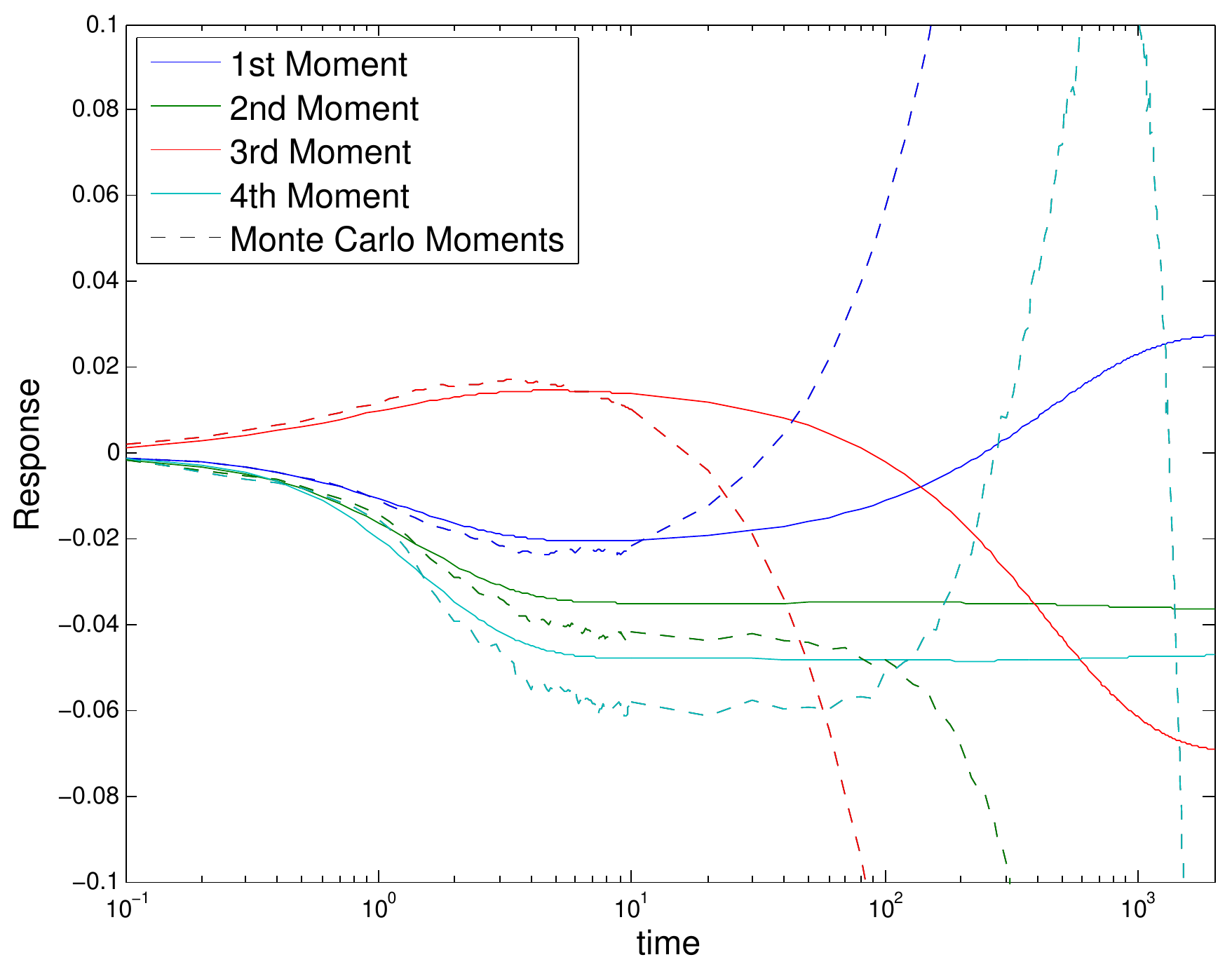}
\includegraphics[width=0.45\textwidth]{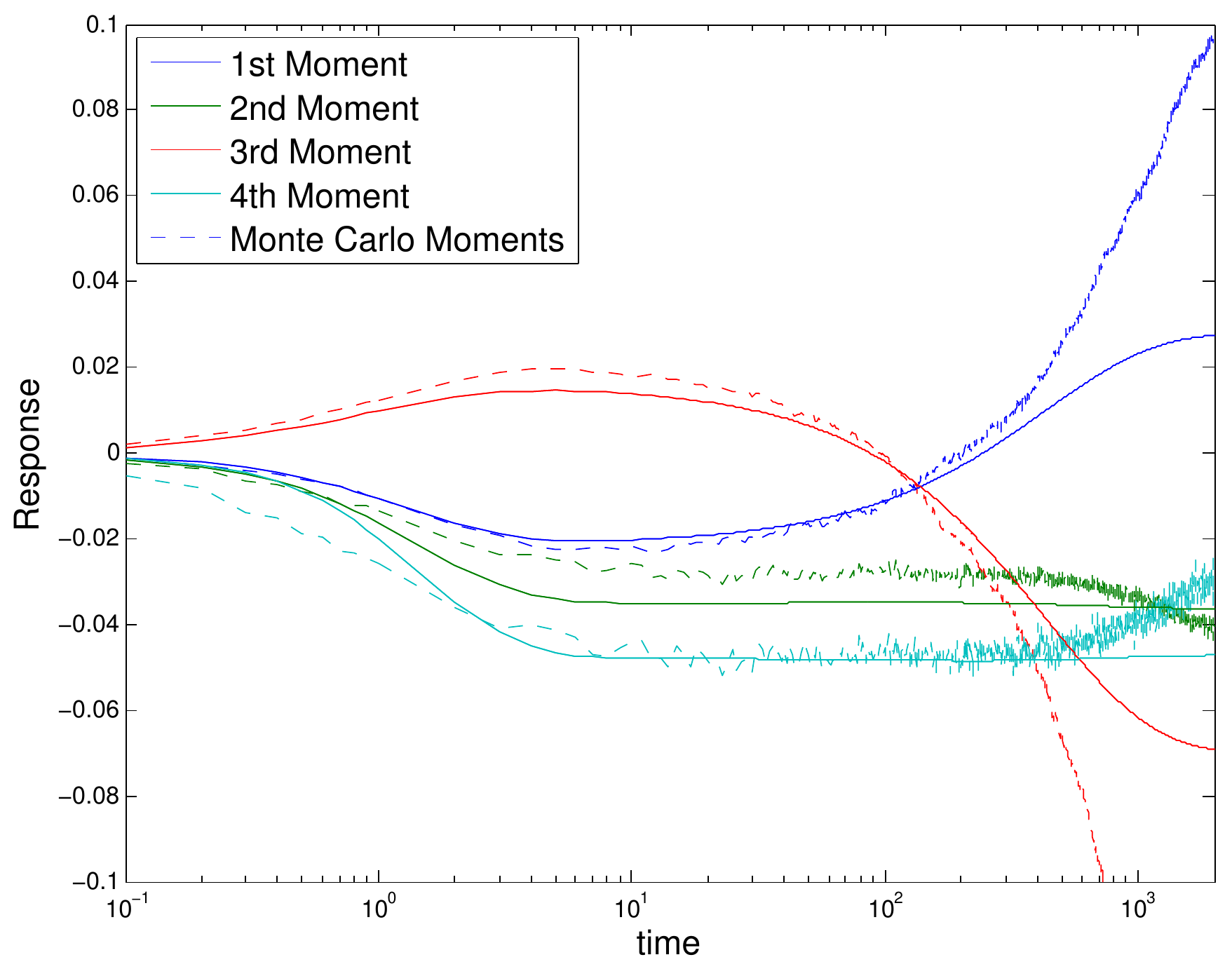}
\caption{\label{doublewellresponsefull} Left: Evolution of the moments compared to the full model with $\epsilon=\sqrt{0.1}$.  Right: Evolution of the moments compared to the full model with $\epsilon=\sqrt{0.005}$. The evolution of the densities is shown in {\tt fig9.mov} for $\epsilon = \sqrt{0.1}$.}
\end{figure}

Finally, we compare the non-parametric response estimate to a Monte-Carlo simulation of the full model \eqref{doublewellfull} with perturbation \eqref{perturbeddoublewell}.  We initialize an ensemble of 100000 members as in the stochastic simulation, and then integrate the perturbed system using a fourth order Runge-Kutta scheme with $\delta t = 0.002$.  As shown in Figure \ref{doublewellresponsefull}, for $\epsilon = \sqrt{0.1}$ this procedure only matches the evolution of the stochastically forced gradient flow for a short time.  The evolution of the densities up to time $t=10$ is shown for $\epsilon=\sqrt{0.1}$ in {\tt fig9.mov}.  By repeating the Monte-Carlo simulation with $\epsilon = \sqrt{0.005}$ and decreasing the integration time step to $\delta t = 0.0005$ we find better agreement for longer time, however the system still does not behave as a gradient flow for very long times as shown in Figure \ref{doublewellresponsefull}.

\section{Summary and Discussion}\label{sec6}

The techniques developed in this paper solve uncertainty quantification problems using only a training data set without any parametric model.  The key tool is the diffusion maps algorithm developed in \cite{diffusion,VB}, which uses the training data to build a non-parametric model for any gradient flow system with isotropic homogeneous stochastic forcing.  The model is represented implicitly through a kernel matrix and which converges to the generator of the dynamical system in the limit of large data \cite{VB}.  Using properties of gradient flow systems, we can use the eigenvectors of this kernel matrix to approximate the eigenfunctions of the corresponding Fokker-Planck operator.  Subsequently, we solved three uncertainty quantification problems, forecasting, filtering, and response estimation, by projecting the associated PDE or SPDE onto this basis of eigenfunctions.  In this basis the system becomes a linear ODE or SDE which we can then truncate onto a finite number of modes and solve analytically or approximate numerically.  

A key advantage of this technique is that all the densities and eigenfunctions are represented by their values on the training data set, which alleviates the need for a grid.  This is important because it is difficult to define a grid on an unknown manifold $\mathcal{M}\subset\mathbb{R}^n$.  While our numerical examples focused on one-dimensional systems on $\mathcal{M}=\mathbb{R}$, the theory of \cite{diffusion,VB} allows the gradient flow system to evolve on any $d$-dimensional manifold of $\mathcal{M}\subset\mathbb{R}^n$.  It is important to note that these techniques provably solve the full nonlinear uncertainty quantification problem, but only in the limit of infinite training data and an infinite number of eigenvectors.  In practice, the amount of data and number of eigenvectors required will depend strongly on the intrinsic dimensionality of the underlying manifold.  For high-dimensional systems the amount of data required will quickly make these techniques computationally infeasible even if such large data sets are available.  

The techniques are also limited by the assumption of a gradient flow system with stochastic forcing which is isotropic and constant over the entire state space.  Concurrent research in \cite{bs:14} suggests that it may be possible to overcome these limitations using different kernels based on non-Euclidean norms which vary over the state space.  Overcoming the restriction to gradient flow systems may also be possible using kernels which are not centered as shown in \cite{bs:14}.  However, one limitation of the results of \cite{bs:14} is the need to estimate the drift and diffusion coefficients that are used to define the appropriate un-centered kernel. Finding provably convergent algorithms for estimating these coefficients can be difficult, particularly on non-compact domains as shown in \cite{VB}.  Instead of explicitly estimating the drift and diffusion coefficients, an alternative approach has been introduced in \cite{BGH14} which implicitly recovers these coefficients. The method of \cite{BGH14} was applied to the problem of forecasting turbulent Fourier modes in \cite{bh:15physd} with small noisy training data sets, whereas here we considered large data sets to verify our approach in the limit of large data.  Crucially, the approach of \cite{BGH14} is still based on representing an initial density in the diffusion coordinates used here, and we believe that the UQ methods developed here can be extended to the more general context of \cite{BGH14,bh:15physd}. Regardless of these developments, the fundamental limitation for nonparametric modeling will likely remain the requirement that the manifold be low-dimensional.


\bibliographystyle{siam.bst}
\bibliography{DMUQbib}


\Appendix
\section{Moments and Nonlinear Response of the Ornstein-Uhlenbeck Process}\label{OUanalysis}

Consider the following linear SDE,
\begin{align}\label{pou} dx = ((\alpha+a)x + b)dt + \sqrt{2D}dW_t, \end{align}
with solutions given by Ornstein-Uhlenbeck processes centered at $\bar x_{\infty} \equiv -b/(\alpha+a)$.  The Fokker-Plank equation for \eqref{pou} is given by,
\begin{align}\label{poufp} \frac{\partial}{\partial t}p^{\delta} = \mathcal{L}^*p^{\delta} &= -\frac{\partial}{\partial x}(((\alpha + a)x+ b)p^{\delta}) + D \frac{\partial^2}{\partial x^2} p^{\delta}, \end{align}
and we can use this equation to find the evolution of the moments.  First, the mean $\bar x(t) = \int x p^{\delta}(x,t)dt$ evolves according to,
\[ \dot{\bar x} = \int x p^{\delta}_t dx = \int x\left(-\frac{\partial}{\partial x}(((\alpha + a)x+ b)p^{\delta}) + D \frac{\partial^2}{\partial x^2} p^{\delta} \right) dx = (\alpha + a) \bar x + b. \]
Using the evolution of the mean, we can then find the evolution of the centered moments $M_n(t) = \int (x-\bar x)^n p^{\delta}(x,t)dx$ using the same strategy,
\begin{align} \dot M_n &= \int (x-\bar x)^n p^{\delta}_t - n \dot{\bar x}(x-\bar x)^{n-1}p^{\delta} \nonumber \\
&= \int n(x-\bar x)^{n-1}((\alpha + a)x+b)p^{\delta} - n(x-\bar x)^{n-1}((\alpha + a)\bar x + b)p^{\delta} + n(n-1)D(x-\bar x)^{n-2}p^{\delta} \nonumber \\
&= n(\alpha + a)M_n + n(n-1)D M_{n-2}.\nonumber
\end{align}
Noting that $M_0 = 1$ and $M_1 = 0$, we have the following solutions for the first four moments,
\begin{align}
\bar x(t) &= (\bar x(0) - \bar x_{\infty})e^{(\alpha + a)t} + \bar x_{\infty} \nonumber \\
M_2(t) &= \left(M_2(0) + \frac{D}{\alpha + a}\right)e^{2(\alpha + a)t} - \frac{D}{\alpha + a} \nonumber \\
M_3(t) &= M_3(0)e^{3(\alpha+a)t} \nonumber \\
M_4(t) &= \left(M_4(0) + \frac{6D}{\alpha + a}\left(M_2(0) + \frac{D}{\alpha+a}\right) - \frac{3D^2}{(\alpha + a)^2} \right)e^{4(\alpha + a)t} \nonumber \\
&\hspace{15pt} - \frac{6D}{\alpha+a}\left(M_2(0) + \frac{D}{\alpha + a}\right)e^{2(\alpha+a)t} + \frac{3D^2}{(\alpha+a)^2}. \nonumber
\end{align}
To find the analytic solutions for the response, $\delta\mathbb{E}[A(x)]$, we take the initial moments from the invariant distribution $p_{\rm{eq}}$ of the unperturbed system.  Thus, $\bar x(0) = 0$, $M_2(0) = -D/\alpha$, $M_3(0) = 0$, and $M_4(0) = 3D^2/\alpha^2$ so plugging in these values and then subtracting them from the moments we have,
\begin{align}
\delta\mathbb{E}[x] &= \bar x_{\infty}\left(1-e^{(\alpha + a)t}\right) =  \frac{-b}{\alpha+a}(1-e^{(\alpha+a)t}) \nonumber \\
\delta\mathbb{E}[(x-\bar x)^2] &= \frac{aD}{\alpha(\alpha + a)}\left(1 - e^{2(\alpha + a)t}\right)\nonumber \\
\delta\mathbb{E}[(x-\bar x)^3] &= 0 \nonumber \\
\delta\mathbb{E}[(x-\bar x)^4] &= \frac{6aD^2}{\alpha(\alpha + a)^2}\left(e^{2(\alpha + a)t} - e^{4(\alpha + a)t}\right) + \left(\frac{3D^2}{\alpha^2}-\frac{3D^2}{(\alpha + a)^2}\right)\left(e^{4(\alpha+a)t} - 1\right). \nonumber
\end{align}

\end{document}